%% file: main.tex
\documentclass[english]{article}
\usepackage[T1]{fontenc}	% encoding czcionek uzytych do wyswietlania (powinien dzialac z UTF)
\usepackage[utf8]{inputenc}	% encoding tekstu zapisanego w pliku (chcyby UTF - dowolne znaki, np. polskie sa ok)
\usepackage{hyperref} 	% ten pakiet sprawia, ze odnosniki w tekscie sa klikalne (latwiej sie pracuje), musi byc przed innymi, inaczej bledy
\usepackage{textcomp}	% ???
\usepackage{amsmath}	% ams - matematyka
\usepackage{graphicx}   % for includegraphics
\usepackage{amssymb}	% ams - matematyka
\usepackage{url}		% urle w cytowaniach
\usepackage{esint}		% ????
\usepackage{babel}		% polskie znaki / polskie dzielenie wyrazow
\usepackage{soul}  		% ten pakiet pozwala uzywac komendy \st{} do przekreslania tekstu
\usepackage{xcolor}  	% for \color
\usepackage{authblk}
\usepackage{tikz}
\usepackage{tikz-cd}
\usepackage{comment}
\usepackage{soul}
\usepackage{enumitem}

\newtheorem{theorem}{Theorem}
\newtheorem{definition}{Definition}
\newtheorem{lemma}[theorem]{Lemma}
\newtheorem{rem}[theorem]{Remark}
\newtheorem{pro}[theorem]{Proposition}

\newtheorem{cor}[theorem]{Corollary}

\newcommand{\grad}{\ensuremath{\nabla}}

\newcommand{\N}{\ensuremath{\mathbb{N}}}
\newcommand{\R}{\ensuremath{\mathbb{R}}}
\newcommand{\B}{\ensuremath{\mathbb{B}}}

\DeclareMathOperator{\Id}{Id}
\newcommand{\cB}{\ensuremath{\overline{\mathbb{B}}}}

\newcommand{\IR}{\mathbb{IR}}
\renewcommand{\O}{\ensuremath{\mathcal{O}}}
\newcommand{\X}{\ensuremath{\mathcal{X}}}
\newcommand{\epsi}{\ensuremath{\varepsilon}}
\DeclareMathOperator{\bd}{\ensuremath{\partial}}
\DeclareMathOperator{\dist}{\ensuremath{dist}}

\DeclareMathOperator{\interior}{int}
\DeclareMathOperator{\domain}{dom}
\DeclareMathOperator{\diam}{diam}
\DeclareMathOperator{\spans}{span}
\newcommand{\cover}[1]{\stackrel{#1}{\Longrightarrow}}

\graphicspath{{./}} % Directory in which figures are stored

\def\qed{{\hfill{\vrule height5pt width5pt depth0pt}\medskip}}

\title{Sharkovskii theorem for infinite dimensional dynamical systems.}

\author{Anna Gierzkiewicz$^{a,b}$}
\author{Robert Szczelina$^{a,c}$}
\affil{$^a$ Institute of Computer Science, Jagiellonian University,\linebreak
ul. \L ojasiewicza 6, 30-348 Krak\'ow, Poland}
\affil{$^b$ \href{mailto:anna.gierzkiewicz@uj.edu.pl}{anna.gierzkiewicz@uj.edu.pl}}
\affil{$^c$ \href{mailto:robert.szczelina@uj.edu.pl}{robert.szczelina@uj.edu.pl, corresponding author}}

\begin{document}

%%%%%%%%%%%%%%%%%%%%%%%%%%%%%%%%%%%%%%%%%%%%%%%%%%%%%%%%%%%%%%%%%%%%%%%%%%%%%%%%
% TITLE/AUTHORS %%%%%%%%%%%%%%%%%%%%%%%%%%%%%%%%%%%%%%%%%%%%%%%%%%%%%%%%%%%%%%%%
%%%%%%%%%%%%%%%%%%%%%%%%%%%%%%%%%%%%%%%%%%%%%%%%%%%%%%%%%%%%%%%%%%%%%%%%%%%%%%%%

\maketitle

\textbf{This work is dedicated to Piotr Zgliczyński on his 60th birthday \vspace{1em}}

%\today

%%%%%%%%%%%%%%%%%%%%%%%%%%%%%%%%%%%%%%%%%%%%%%%%%%%%%%%%%%%%%%%%%%%%%%%%%%%%%%%%
% ABSTRACT %%%%%%%%%%%%%%%%%%%%%%%%%%%%%%%%%%%%%%%%%%%%%%%%%%%%%%%%%%%%%%%%%%%%%
%%%%%%%%%%%%%%%%%%%%%%%%%%%%%%%%%%%%%%%%%%%%%%%%%%%%%%%%%%%%%%%%%%%%%%%%%%%%%%%%
\begin{abstract}
We present an adaptation of a relatively simple topological
argument to show the existence of many periodic orbits in an
infinite dimensional dynamical system, provided that
the system is close to a one-dimensional map 
in a certain sense.
Namely, we prove a
Sharkovskii-type theorem: if the system has a periodic
orbit of basic period $m$, then it must have all periodic orbits
of periods $n \triangleright m$, for $n$ preceding
$m$ in Sharkovskii ordering. The assumptions of the theorem
can be verified with computer assistance, and we 
demonstrate the application of such an argument 
in the case of Delay Differential Equations (DDEs):
we consider the R\"ossler ODE system
perturbed by a delayed term and we show that it retains periodic orbits of
all natural periods for fixed values of parameters. 
\\[2mm]
\textbf{Keywords:} Sharkovskii theorem, computer-assisted proof, periodic orbits, 
delay differential equations (DDEs)
% \\[2mm]
% \textbf{2010 Mathematics Subject Classification:} 
% % 34C25, % periodic orbits in ODEs
% 34K13, % periodic orbits to FDEs
% 34K23, % Complex (chaotic) behavior of solutions to functional-differential equations 
% 34K27, % Perturbations of functional-differential equations
% 37C25, % Fixed points and periodic points of dynamical systems; fixed-point index theory; local dynamics
% 37B10, % Symbolic dynamics
% 65G20, % Algorithms with automatic result verification
\end{abstract}

%%%%%%%%%%%%%%%%%%%%%%%%%%%%%%%%%%%%%%%%%%%%%%%%%%%%%%%%%%%%%%%%%%%%%%%%%%%%%%%%
% ALL SECTIONS %%%%%%%%%%%%%%%%%%%%%%%%%%%%%%%%%%%%%%%%%%%%%%%%%%%%%%%%%%%%%%%%%
%%%%%%%%%%%%%%%%%%%%%%%%%%%%%%%%%%%%%%%%%%%%%%%%%%%%%%%%%%%%%%%%%%%%%%%%%%%%%%%%
\input{1_intro.tex}

\input{2_covrel_sharkovskii}

\input{3_perturbation}

\input{4_computer_assisted}

\input{5_conclusions}

%%%%%%%%%%%%%%%%%%%%%%%%%%%%%%%%%%%%%%%%%%%%%%%%%%%%%%%%%%%%%%%%%%%%%%%%%%%%%%%%
% ACKNOWLEDGEMENT %%%%%%%%%%%%%%%%%%%%%%%%%%%%%%%%%%%%%%%%%%%%%%%%%%%%%%%%%%%%%%
%%%%%%%%%%%%%%%%%%%%%%%%%%%%%%%%%%%%%%%%%%%%%%%%%%%%%%%%%%%%%%%%%%%%%%%%%%%%%%%%

\section*{CRediT authorship contribution statement}

\noindent {\bf Anna Gierzkiewicz}: Conceptualisation, Methodology, Software, Formal analysis, Writing -- original draft, Writing -- review \& editing. 

\noindent {\bf Robert Szczelina}: Conceptualisation, Methodology, Software, Formal analysis, Writing -- original draft, Writing -- review \& editing, Funding acquisition, Project administration.

\section*{Declaration of Competing Interest}

The authors declare that they have no known competing financial 
interests or personal relationships that could have appeared 
to influence the work reported in this paper.

\section*{Acknowledgements}
AG and RS research was funded in whole by National Science Centre, Poland, 2023/49/B/ST6/02801. 
For the purpose of Open Access, the author has applied a CC-BY public copyright licence 
to any Author Accepted Manuscript (AAM) version arising from this submission.

\section*{Dedication}
The authors would like to thank Piotr Zgliczyński, our supervisor, teacher and great friend, for his constant support and motivation that guided us in all our scientific endeavours. 

%%%%%%%%%%%%%%%%%%%%%%%%%%%%%%%%%%%%%%%%%%%%%%%%%%%%%%%%%%%%%%%%%%%%%%%%%%%%%%%%
% BIBLIOGRAPHY %%%%%%%%%%%%%%%%%%%%%%%%%%%%%%%%%%%%%%%%%%%%%%%%%%%%%%%%%%%%%%%%%
%%%%%%%%%%%%%%%%%%%%%%%%%%%%%%%%%%%%%%%%%%%%%%%%%%%%%%%%%%%%%%%%%%%%%%%%%%%%%%%%
%\bibliographystyle{plain} 	% in alphabetical order (firs author?)
\bibliographystyle{unsrt} 	% in citation order % THIS IS REQUIRED BY CNSNS
\bibliography{dde}

\end{document}

%% file: 1_intro.tex
%%%%%%%%%%%%%%%%%%%%%%%%%%%%%%%%%%%%%%%%%%%%%%%%%%%%%%%%%%%%%%%%%%%%%%%%%%%%%%%%
% SECTION %%%%%%%%%%%%%%%%%%%%%%%%%%%%%%%%%%%%%%%%%%%%%%%%%%%%%%%%%%%%%%%%%%%%%%
%%%%%%%%%%%%%%%%%%%%%%%%%%%%%%%%%%%%%%%%%%%%%%%%%%%%%%%%%%%%%%%%%%%%%%%%%%%%%%%%
\section{Introduction}

Our motivation is the recent paper \cite{GZ2022} where the authors show a 
generalisation of Sharkovskii Theorem for multidimensional perturbations 
of a one-dimensional map having an attracting periodic orbit. 
As an application, they study the R\"ossler system \cite{Rossler76} 
\begin{equation}\label{eq:rossler}\tag{R\"o}
\begin{pmatrix}
x' \\ y' \\ z' 
\end{pmatrix}
= f_{a,b}(x,y,z):=
\begin{pmatrix}
-y-z,
\\
b y+x,
\\
z (x-a)+b,
\end{pmatrix}
\end{equation}
or, more precisely, a (local) Poincar\'e map $P : S \to S$
for \eqref{eq:rossler} on the section $S$
defined by $x = 0, x' > 0$.
For some fixed values of parameters $a$, $b$ 
the attractor of $P$
is a 3-, 5- or 6-periodic point $x$ for $P$. 
The generalisation (Th.\ \ref{thm:sh-manyD}) of Sharkovskii theorem states that 
$P$ has also periodic points
of all periods greater in the sense of Sharkovskii ordering 
(see Section~\ref{sec:sharkovskii}) than the period of $x$. 
In view of the above, informally, 
we will call such a situation the \emph{Sharkovskii property}
of $P$. The exact definition will be given later with the
statement of Theorems~\ref{thm:sh-manyD}~and~\ref{thm:sh-infD} (Def.~\ref{def:SharProp}).
Of course, each periodic point of $P$ corresponds to a periodic orbit for \eqref{eq:rossler},
so in that sense we will be writing that the R\"ossler
system \eqref{eq:rossler} has the Sharkovskii property too. We will
call the periodic orbit $x$ for  \eqref{eq:rossler} $m$-periodic
if the corresponding initial value $x(0) \in S$ is an $m$-periodic point for $P$.
The aim of the present paper is to show that a method from  \cite{GZ2022}
combined with the notion of covering relations from \cite{szczelina-zgliczynski-focm-2}
leads to a relatively easy extension to DDEs and other infinite-dimensional 
dynamical systems in the form of Theorem~\ref{thm:sh-infD}, which is a 
generalisation of \cite[Th.\ 2]{GZ2022} for compact mappings on Banach spaces. 

As an application, we apply Theorem~\ref{thm:sh-infD} in two ways. 
First, we show that a sufficiently small perturbation of an ODE to
a DDE preserves Sharkovskii property. Secondly,  in Theorem~\ref{thm:main-numerical}, 
we apply rigorous numerical methods of \cite{szczelina-zgliczynski-focm-2} 
to prove that the system \eqref{eq:rossler} with explicitly given $f_{5.25, 0.2}$ 
and right hand side perturbed by the delayed term $10^{-4} f_{5.25, 0.2}(t-0.5)$ 
also has the Sharkovskii property. 
This proof is computer-assisted \cite{dde-ros-sha-codes}, where 
by a computer-assisted proof we understand a rigorous mathematical 
reasoning with some difficult to calculate by hand part replaced by a 
potentially huge number of simpler conditions that are easy to verify 
by a computer. In our computations we use the \texttt{CAPD} library 
for \texttt{C++} programming language \cite{capd-article} and
its extension to handle general DDEs \cite{szczelina-zgliczynski-focm-2}.

We would like to stress the fact that the method presented in 
this paper is rather general and does not technically 
limit to perturbations of DDEs. It can be applied to any setting
where the dynamical system might be used to generate some
form of \emph{covering relations}, e.g. in PDEs as was shown in 
\cite{wilczak-pde-1,zgliczynski-pde-3}. Of course applicability 
might be limited due to the costly computations 
(\emph{curse of dimensionality}, see discussion 
in Section~\ref{sec:conclusions}). 

There are some works that tackle problem of persistence of 
periodic orbits in ODEs perturbed with delayed term, e.g. 
\cite{delallave-sdde-ode}, even in case of state-dependent delays. 
Also, the fact that \emph{small delays do not matter} for the 
dynamics is well known quantitatively 
\cite{dde-perturb-kurzweil1971,wojcik-zgliczynski-dde-chaos}. 
Our work differ in that we show the persistence of the Sharkovskii 
property -- existence of infinitely many periodic orbits, 
and the persistence of qualitative description of the dynamics.

Some papers use Sharkovskii theorem to investigate dynamics 
of DDEs (see \cite{Ivanov1992} and references therein),  
but it seems that those works mainly deal with  
perturbations of 1-dimensional interval maps for which the Sharkovskii 
property holds. In contrast to this approach, with our methods 
we can start with the infinite-dimensional system
and verify conditions that guarantee the dynamics 
is close to a 1 dimensional map.

In recent years, there appeared some new computer assisted protocols 
that can lead to a proof of complicated dynamics in general DDEs 
such as the famous Mackey-Glass equation 
\cite{numerical-homoclinic-mg, complicated-periodic-mg-limit}. 
We hope that our work may enable another possible approach to solve 
this open problem. 

% SUBSECTION %%%%%%%%%%%%%%%%%%%%%%%%%%%%%%%%%%%%%%%%%%%%%%%%%%%%%%%%%%%%%%%%%%%%%%%%
\subsection*{Outline of the paper}
In Section~\ref{sec:sharkovskii} we present the necessary theory and
the main theorem of this paper (Theorem \ref{thm:sh-infD}). It is worth 
noting that this Theorem is quite general and might be applied 
in other infinite-dimensional scenarios, not only DDEs. 
In Section~\ref{sec:perturbation} we provide Theorem \ref{thm:perturbation}
to show that Delay Differential Equations being a sufficiently small 
perturbation of an ODE model preserves the Sharkovskii property (Def.~\ref{def:SharProp}). 
The theory that `small delays don't matter'
is well known \cite{dde-perturb-kurzweil1967,dde-perturb-kurzweil1971}, 
but we wanted to provide more details, including some estimates of
how small the delay needs to be for the theorems to hold. 
This was done to compare with the results of the Section~\ref{sec:CAP},
where we show the existence of complicated dynamic in a 
R\"ossler system perturbed with a delayed term for fixed and
reasonably macroscopic values of parameters. In Section~\ref{sec:conclusions}
we give some possible future applications of the method. 

% SUBSECTION %%%%%%%%%%%%%%%%%%%%%%%%%%%%%%%%%%%%%%%%%%%%%%%%%%%%%%%%%%%%%%%%%%%%%%%%
\subsection*{Glossary and arrangements}
\begin{itemize}
	\item The closure, the interior, and the boundary of a topological set 
	$A\subset\mathbb{R}^N$ are denoted by $\overline{A}$, $\interior A$, 
	and $\partial A$, respectively.
	\item $\B_N$, $\cB_N$ are the unit open and closed balls in $\R^N$ with 
	the maximum norm $\|\cdot\|=\|\cdot\|_{\max}$. In general, we will also 
	use $\B_N(a,R)$, $\cB_N(a,R)$ for balls of radius $R$ centred in $a$ 
	and $\B_N(R) = \B_N(0,R)$, $\cB_N(R) = \cB_N(0,R)$. If $\mathcal{X}$ 
	is a Banach space, then by $\B_{\mathcal{X}}(a, R)$, etc. we denote 
	the respective ball in that space. Finally, for any set $A$ in some 
	Banach space $\mathcal{X}$, we denote by $\B(U, r)$, the open neighbourhood 
	of $U$ of radius $r$, i.e. $\B(U, r) = \{u \in U:  \dist(u, U) < r \}$.
	\item For a vector $v$ in a Hilbert  space $\mathcal{X}$ by 
	$\spans(v_1, \ldots, v_k)$ we denote the subspace spanned by vectors 
	$v_i$. If $U$ is a subspace of $\mathcal{X}$ then by $U^\bot$ we 
	denote the maximal subspace in $\mathcal{X}$ that is orthogonal 
	to all vectors in $U$ w.r.t. the scalar product in $\mathcal{X}$.
	\item The projection operator onto the $i$th coordinate in $\mathbb{R}^N$ 
	is denoted by $\pi_i$. By $\pi_{\X_i}$ we denote the projection of a space 
	on its subspace $\X_i$. By $\pi_{i_1,i_2,\ldots,i_k}$ we denote the 
	projection $\R^d \to \R^k$ given by $(\pi_{i_1}, \ldots, \pi_{i_k})$, 
	and by $\pi_{j \ge k}$ we denote the projection $\pi_{k, k+1, \ldots, d}$. 
	The dimension $d$ will be known from the context.
	\item By an `$n$-periodic orbit' or `$n$-periodic point', we understand 
	an orbit or a point with fundamental period $n$.
	\item Whenever the natural counter is limited, \textit{eg.} 
	$i\in\{0,1,\dots,n-1\}$, then by `$i$' we mean `$i\pmod{n}$'. 
	In particular, if $i=n-1$, then $J_{i+1}=J_0$.
\end{itemize}

%% file: 2_covrel_sharkovskii.tex
%%%%%%%%%%%%%%%%%%%%%%%%%%%%%%%%%%%%%%%%%%%%%%%%%%%%%%%%%%%%%%%%%%%%%%%%%%%%%%%%
% SECTION %%%%%%%%%%%%%%%%%%%%%%%%%%%%%%%%%%%%%%%%%%%%%%%%%%%%%%%%%%%%%%%%%%%%%%
%%%%%%%%%%%%%%%%%%%%%%%%%%%%%%%%%%%%%%%%%%%%%%%%%%%%%%%%%%%%%%%%%%%%%%%%%%%%%%%%
\section{\label{sec:sharkovskii}Sharkovskii-type theorem in infinite dimensions}

In this chapter we will show how to extend the Sharkovskii
theorem to compact maps in infinite dimensional spaces,
which are sufficiently close to one-dimensional expanding maps. 
For completeness we will also remind the basic notion of
Sharkovskii ordering, the Theorem itself and the sketch
of the proof in the finite dimensional case. 

\subsection{Sharkovskii Theorem}

\begin{theorem}[Sharkovskii, \cite{ShU}]\label{th:shar}
Define an ordering `$\triangleleft$' of natural numbers:
\begin{equation}\label{eq:order}
\begin{array}{r@{\ \triangleleft\ }c@{\ \triangleleft\ }c@{\ \triangleleft\ }c@{\ \triangleleft\ }l}
3\triangleleft 5 \triangleleft 7 \triangleleft 9 \triangleleft\ \dots & 2\cdot 3 & 2 \cdot 5 & 2 \cdot 7 & \dots
\\
 \dots & 2^2\cdot 3 & 2^2 \cdot 5 & 2^2 \cdot 7 & \dots
\\
 \dots & 2^3\cdot 3 & 2^3 \cdot 5 & 2^3 \cdot 7 & \dots
\\
 \multicolumn{5}{c}{\dots \dots \dots \dots \dots \dots \dots \dots \dots \dots }
\\
 \dots & 2^{k+1} & 2^k & 2^{k-1} & \ldots\ \triangleleft 2^2 \triangleleft 2 \triangleleft 1.
\\
\end{array}
\end{equation}
Let $f: I \to \mathbb{R}$ be a continuous map of an interval. If $f$ has an $n$-periodic point and $n\triangleleft m$, then $f$ also has an $m$-periodic point.
\end{theorem}
In particular, in the case $n=3$, if the map has a $3$-periodic orbit, then it also has orbits of all natural periods (\emph{Period 3 implies chaos}, \cite{LiYorke}). 

The Sharkovskii Theorem \cite{ShU} describes a situation which may occur for one dimensional maps exclusively, and is often considered a chaotic-like behaviour of a map \cite{LiYorke}. Maybe this unusual, fascinating ordering of naturals is what causes the unabated interest in this field and plenty of different proofs of the Theorem that appeared up to date. One of the simplest and most intuitive proofs was given by Burns and Hasselblatt \cite{burns}, and was, in addition, a great help in proving an extension of the Theorem for multidimensional perturbations of one-dimensional maps \cite{GZ2022}.
Let us recall shortly the basic tools used there.

The following Definition of covering relation between intervals comes from \cite[Def.\ 9]{GZ2022} and is a slight modification of the one that is often used in Sharkovskii Theorem proving, see \textit{e.g.}\ Definition \cite[Def.\ 2.1]{burns}. In this paper, we will omit the adjective `proper'.

\begin{definition}[(Proper) covering relation between intervals]\label{def:Oforced_covering}
Let $I, J \subset \mathbb{R}$ be closed intervals and consider a continuous map $f:I\to\mathbb{R}$.
We say that the interval $I$ \emph{covers} the interval $J$ and denote it by $I \overset{f}{\rightarrow} J $ or simply  $I \rightarrow J $, if 
\begin{equation}\label{eq:proper_cover}
\min f(\partial I) \leq \min J \text{ \qquad and \qquad } \max J \leq \max f(\partial I).
\end{equation}
\end{definition}

The above relation fulfils the next Itinerary Lemma, which is crucial in proving the existence of periodic orbits in the proofs of Sharkovskii's theorem:
\begin{theorem}[Itinerary Lemma]
\label{th:1d-covering}
Let $f: I \to \mathbb{R}$ be a continuous map on an interval $I \subset\mathbb{R}$. Assume that we have a sequence of intervals $J_j \subset I$ for $j=0,\dots,m-1$ such that
\begin{equation}\label{eq:1d_loop}
  J_0  \overset{f}{\longrightarrow} J_1 \overset{f}{\longrightarrow} J_2 \overset{f}{\longrightarrow} \dots \overset{f}{\longrightarrow} J_{m-1} \overset{f}{\longrightarrow} J_0.
\end{equation}
Then there exists a point $x \in J_0$, such that $f^j(x) \in J_j$ for $j=1,\dots,m-1$ and $f^m(x)=x$.
\end{theorem}
A point which fulfils the thesis of Theorem \ref{th:1d-covering} is said to \emph{follow the loop} \eqref{eq:1d_loop}. A loop of length $m$ (such as \eqref{eq:1d_loop}), we will call shortly an $m$-loop.

\begin{definition}[\cite{burns}, Lemma 2.6]\label{def:basic_loop}
We call a loop \eqref{eq:1d_loop} a \emph{non-repeating loop}, if it fulfils the following conditions:
\begin{enumerate}
	\item it is not followed by any endpoint $x\in\bigcup_{i=0}^{m-1} \partial J_i$;
	\item $\displaystyle \interior J_0 \cap \bigcup_{i=1}^{m-1} J_{i} = \varnothing$.
\end{enumerate}
\end{definition}

\begin{rem}[\cite{burns}, Lemma 2.6]\label{lem:non-repeating}
If a loop \eqref{eq:1d_loop} is non-repeating, then every point following it has least period $m$.
\end{rem}

Let now $f: I \to \mathbb{R}$ be a continuous map on an interval $I \subset\mathbb{R}$ and $\O\subset I$ be an $n$-periodic orbit for $f$. 
By an \emph{$\mathcal{O}$-interval} we understand any interval $J\subset [\min\O, \max\O]$ of positive length with the endpoints in $\mathcal{O}$, that is, $\partial J \subset \mathcal{O}$.

%\subsubsection{The Proposition}
The following theorem is Lemma 11 from \cite{GZ2022} and is a slight modification of Proposition 6.1 of \cite{burns} for proper covering relations.
\begin{theorem}[\cite{GZ2022}]\label{th:Burns}
	For every number $m$ succeeding $n$ in Sharkovskii's order \eqref{eq:order} (\textit{i.e.} $n\triangleleft m$) there exists a non-repeating %$\O$-forced 
	$m$-loop of $\mathcal{O}$-intervals $J_i$, $i=0,\dots, m-1$:
	\begin{equation}
	  J_0  \overset{f}{\longrightarrow} J_1 \overset{f}{\longrightarrow} J_2 \overset{f}{\longrightarrow} \dots \overset{f}{\longrightarrow} J_{m-1} \overset{f}{\longrightarrow} J_0\text{,}
	\end{equation}
	which proves the existence of an $m$-periodic point for $f$ in the interval $ [\min \mathcal{O},\max\mathcal{O}]$.
\end{theorem}

\subsection{Covering relations}

%Gidea, Zgliczy\'nski: \cite{GZ,Z}

%Miranda: \cite{Mir}

 In order to have an analogue of the Itinerary Lemma in many dimensions we need a good notion of covering where we have a direction of possible expansion
and an apparent contraction in other directions.

For this end  we recall  the notion of covering for h-sets in $\mathbb{R}^N$ from \cite{PZszarI,PZmulti,GZ2021}.
It is, in fact, a particular case of the similar notion from \cite{GZ}, but with exactly one exit (or 'unstable')  direction. For such a covering relation a special version of Itinerary Lemma is true (Th.\ \ref{thm:periodic-covering-FINITE}  below) and we use it to prove the existence of periodic points for multidimensional maps.

\begin{definition}[Definition~1 in \cite{GZ}, for the case $u_{N}=1$]
\emph{A (horizontal) h-set} $N$ in $\R^{d_N}$, $d_N\geq 2$ is a pair $(|N|,c_N)$ where:
\begin{itemize}
\item $|N|$ is a compact subset of $\R^{d_N}$ (called a \emph{support} of the h-set);
\item a homeomorphism $c_N : \R^{d_N} \to \R^{d_N} = \R \times \R^{d_N-1}$ is such that
\begin{equation*}
c_N(|N|) = [-1,1] \times \cB_{d_N-1}.
\end{equation*}
\end{itemize}
We define also
\begin{eqnarray*}
N_c & = & [-1,1] \times \cB_{d_N-1}  \\
N^-_c & = & \{-1,1\} \times \cB_{d_N-1} \\
N^+_c & = & [-1,1] \times \bd  \cB_{d_N-1} \\
N^- & = & c_N^{-1}(N^-_c)  \text{\quad (the exit set)}\\
N^+ & = & c_N^{-1}(N^+_c) \text{\quad  (the entry set).}
\end{eqnarray*}
\end{definition}
Basically, here an h-set $N$ is a product of a closed interval and a closed $(d_N-1)$-dimensional ball
in an appropriate coordinate system. Note that the exit set has always two connected components homeomorphic to $\cB_{d_N-1}$. One can just simply imagine a hyper-cuboid with its left and right sides distinguished.
The number $d_N-1$
is the entry dimension (nominally stable), and the exit dimension will be equal to one in this paper. Of course, one can set the definition also for any natural number of unstable dimensions (see \cite{GZ}), but here we skip it on purpose to simplify the next definitions. We will also usually identify the support
$|N|$ with the h-set $N$ (e.g. we will write $f(N)$ instead of $f(|N|)$.

Now we define a horizontal covering relation between h-sets through a continuous map $f$.

\begin{definition}[cf. Definition~2 in \cite{GZ}]
\label{def:singlevalued-covering-FINITE}
Assume $N$, $M$ are two h-sets. Let $f : |N| \to \R^{d_M}$ be a continuous map.
We say that $N$ $f$-covers $M$, denoted by:
\begin{equation*}
N \cover{f} M
\end{equation*}
iff there exists continuous homotopy $H : [0,1] \times |N| \to R^{d_M}$
satisfying the following conditions:
\begin{itemize}
\item $H(0, \cdot) = f$;
\item $H([0,1], N^-) \cap M = \varnothing$;
\item $H([0,1], N) \cap M^+ = \varnothing$;
\item there exists %a linear map $A : \R \to \R$ 
 $a\in\mathbb{R}\setminus[-1,1]$ such that
\begin{eqnarray*}
H_c(1, (p,q)) & = & (a p, 0) 
%\\
%A(\{-1,1\}) & \subset & \R \setminus [-1,1]
\end{eqnarray*}
where $H_c(t, \cdot) = c_M \circ H(t, \cdot) \circ c_N^{-1}$ is
the homotopy expressed in `basic' coordinates.
\end{itemize}
\end{definition}

\begin{figure}[ht]
\begin{center}
	\includegraphics[width=0.5\textwidth]{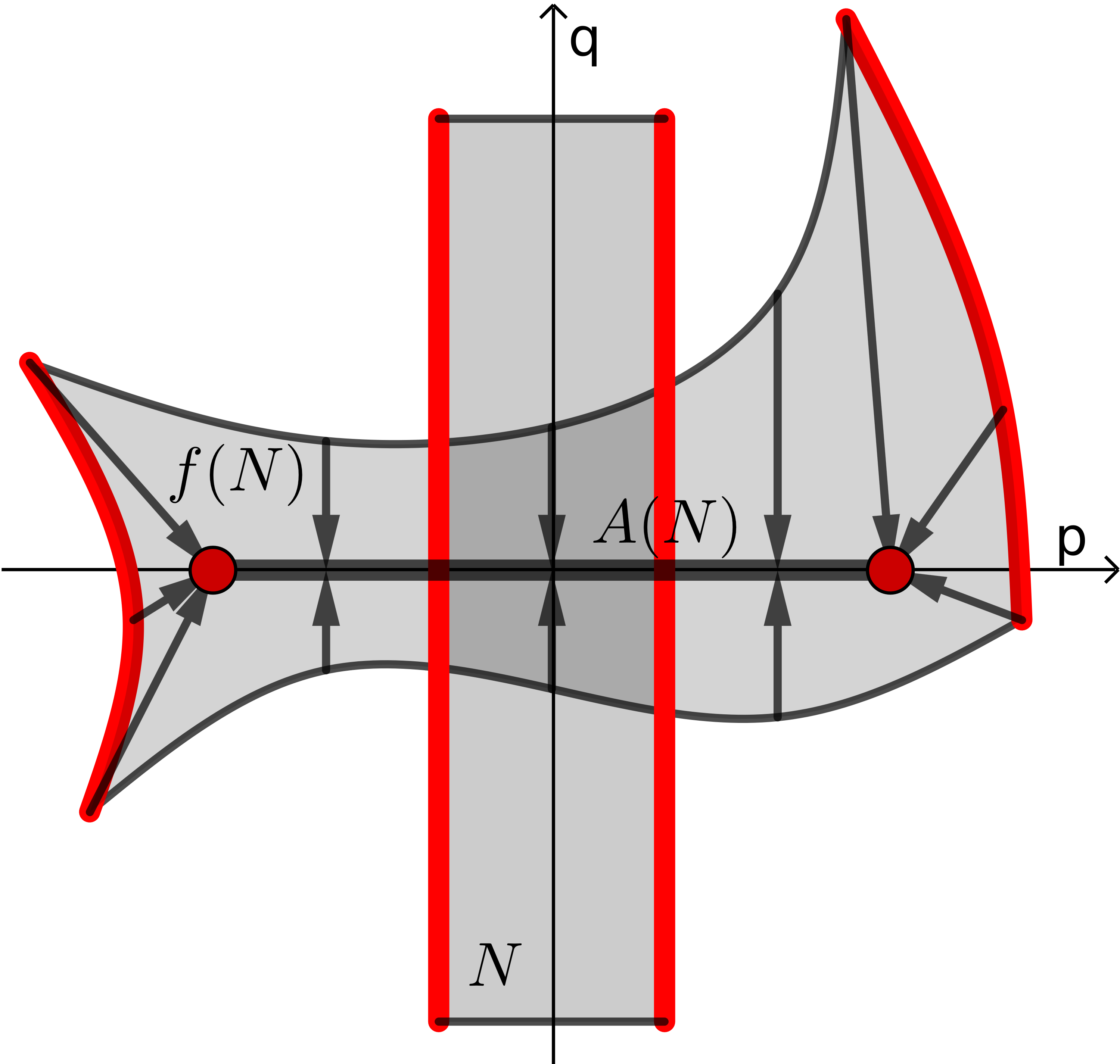}
	\caption{\label{fig:cover}Horizontal self-covering $N \overset{f}{\Longrightarrow} N$. The image $f(N)$ is homotopically equivalent to the image $A(N)$ through a linear map $A(p,q)=(ap,0)$.%
	}
\end{center}
\end{figure}

\begin{figure}[ht]
\begin{center}
	\includegraphics[width=0.9\textwidth]{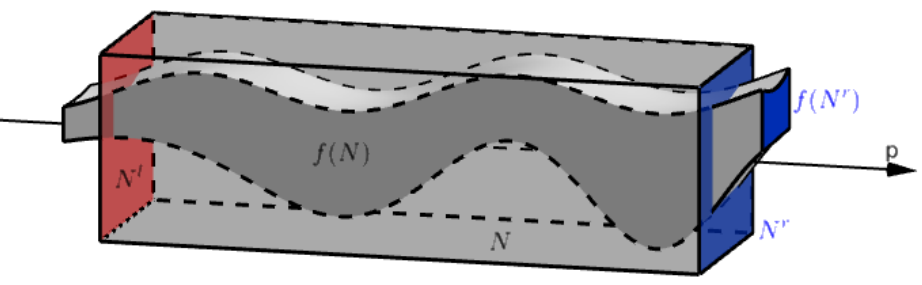}
	\caption{\label{fig:cover3d}Horizontal self-covering $N \overset{f}{\Longrightarrow} N$ in $\R^3$. Note that the `right edge' $N^r$ (blue) is mapped to the right of the h-set and similarly the left edge $N^l$ (red) is mapped to the left.}
\end{center}
\end{figure}

See Fig.\ \ref{fig:cover} for the illustration of horizontal covering in $\mathbb{R}^2$ and Fig.\ \ref{fig:cover3d} for covering in $\mathbb{R}^3$.

\begin{rem}
Note that in this simple case of one unstable direction, the covering $N \cover{f} M$ means that the image of $N$ is stretched along the unstable direction of $M$, contracted in all the other dimensions, and the left and right sides of $N$ must be mapped to the left and to the right of the left and right sides of $M$, respectively (possibly, in the reversed order).
What is important for computer-assisted proofs, the above conditions can be easily checked with the use of computer via interval arithmetic and `$<$', `$>$' relations.
\end{rem}

The following theorem is a version of Itinerary Lemma (Theorem~\ref{th:1d-covering}) for the horizontal covering relation.

\begin{theorem}[\cite{PZszarI}, more generally \cite{GZ}]
\label{thm:periodic-covering-FINITE} 
Suppose that there occurs a loop of $m$ horizontal $f$-coverings between h-sets $M_i\subset\mathbb{R}^d$, $i=0,\dots,m-1$:
\begin{equation}\label{eq:loop}
M_0 \overset{f}{\Longrightarrow} M_1\overset{f}{\Longrightarrow} \dots \overset{f}{\Longrightarrow} M_{m-1} \overset{f}{\Longrightarrow} M_m = M_0\text{,}
\end{equation}
then there exists $x\in M_0$ such that $f^m(x)=x$ and
\[
\text{for } i=1,\dots ,m-1: \qquad f^i(x)\in M_i.
\]
\end{theorem}

% SUBSECTION %%%%%%%%%%%%%%%%%%%%%%%%%%%%%%%%%%%%%%%%%%%%%%%%%%%%%%%%%%%%%%%%%%%
\subsection{The case of infinite dimensional space}

The proof of the Itinerary Lemma \ref{thm:periodic-covering-FINITE} for finite-dimensional covering relations and its various versions \cite{Z,GZ} relies on Miranda Theorem \cite{Mir} or, more generally, Brouwer Fixed Point Theorem, using the notion of fixed point index. In order to define infinite-dimensional h-sets analogues and covering relations between them one needs a similar construct, for instance the Leray--Schauder degree \cite{Schauder}. Then, it is possible to generalise this tool to maps on compact sets, which has applications to PDEs \cite{chaos-kuramoto}, or to compact maps as in \cite{szczelina-zgliczynski-focm-2}, which is of use in our case of Delay Differential Equations, or possibly in a more general setting of Functional Differential Equations.

Here, we follow the definitions from \cite{szczelina-zgliczynski-focm-2}: we will also work on spaces of the type
$\X = \X_1 \oplus \X_2$, where
$\X_1$ is finite-dimensional ($\X_1 \equiv \R^M$) and
$\X_2$ infinite-dimensional. In our application to DDEs,
$\X = C^n_p = \R^{M(d,p,n)} \times (C^0([0,h], \R))^{d \cdot p}$.
The notion of those spaces will be introduced later in Section~\ref{sec:CAP}.

First, recall the notion of \emph{h-set with tail} in the case of one expanding direction (sometimes referred to as \emph{horizontal h-set}).

\begin{definition}[\cite{szczelina-zgliczynski-focm-2}]
\label{def:h-set-with-tail}
Let $\X = \R^{d_N}\times \X_2$ be a real Banach space.

A (horizontal) \emph{h-set with tail} is a pair $N = (N_1, |N_2|)$ where
\begin{itemize}
\item
	$N_1$ is an h-set in $\R^{d_N}$ (the \emph{head}),
\item
	$|N_2| \subset \X_2$ is a closed, convex and bounded set (the \emph{tail}).
\end{itemize}
Additionally, we set %$u_N = u_{N_1}$, 
$|N| = |N_1| \times |N_2|$,
$c_N = (c_{N_1}, \Id)$ and
\begin{eqnarray*}
N_c & = & c_N\left(|N|\right) \ = \ N_{1,c} \times |N_2| \ = \\
    & = & [-1,1] \times \cB_{{d_N}-1}\times |N_2|.
\end{eqnarray*}
\end{definition}
We will just say shortly that $N$ is an h-set when context is clear. In some sense, a finite-dimensional h-set $N$ in $\R^{d_N}$ can be viewed as an h-set with trivial tail $\R^{0} = \{ 0 \}$.

Now, recall the generalisation of covering relation for h-sets with tails.
In what follows, by a \emph{compact mapping} on a Banach space $\mathcal{X}$ 
we understand a map $f : \domain(f) \subset \mathcal{X} \to \mathcal{X}$ such
that for any bounded $Y \subset \domain(f)$, the closure of its image $\overline{f(Y)}$ is compact in $\mathcal{X}$.

\begin{figure}[ht]
\includegraphics[width=\textwidth]{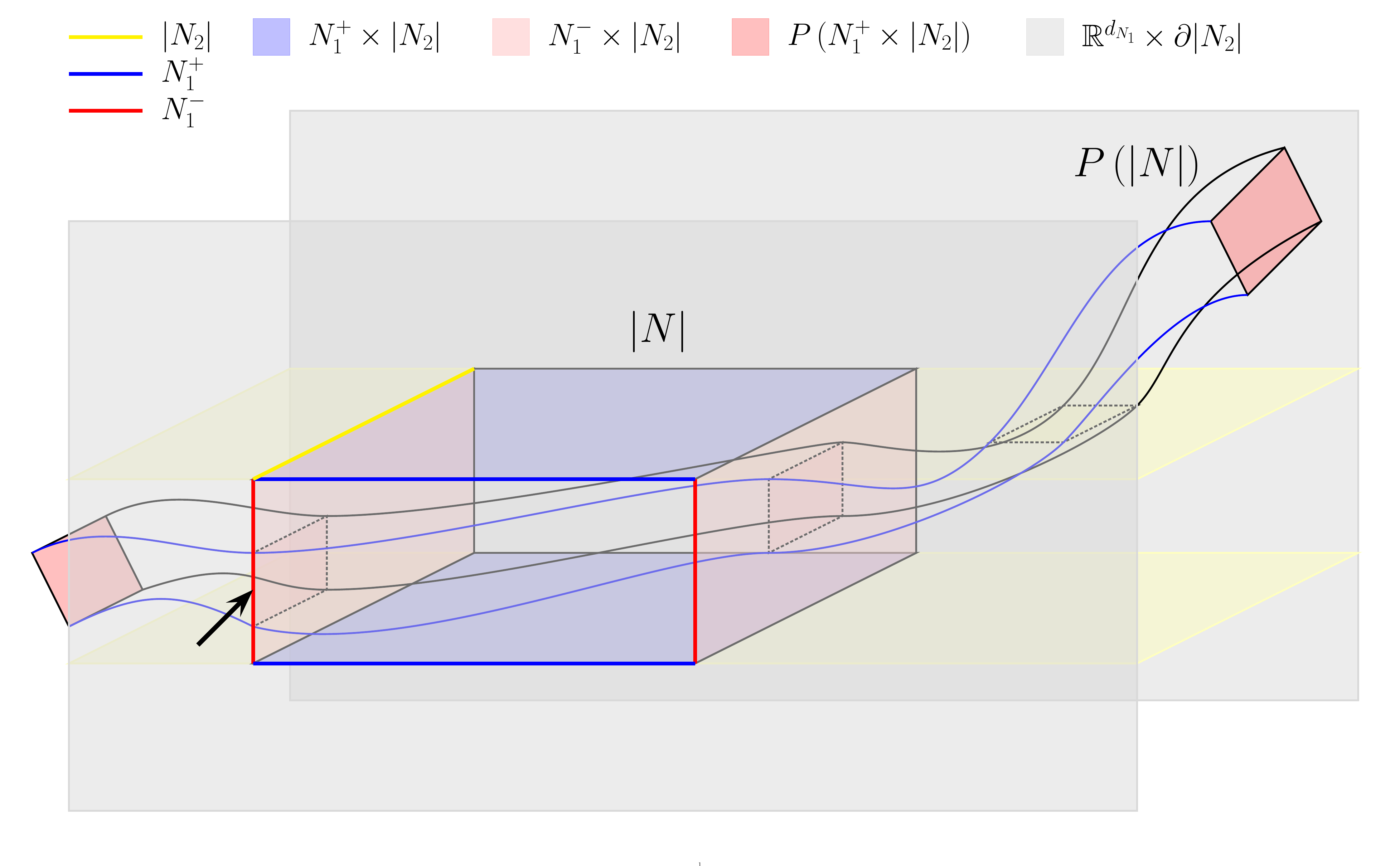}
\caption{\label{fig:covrel}An example of a self-covering relation $N \cover{P} N$
on an h-set with tail $N = (N_1, |N_2|)$, $d_N = 2$. The head $N_1$ is the rectangle in the front, with red and blue sides (the exit and entrance set, respectively). The tail $|N_2|$
is closed and convex in a potentially infinite dimensional space (yellow dimension, the `depth' of the picture). \\
}
\end{figure}

\begin{definition}[\cite{szczelina-zgliczynski-focm-2}]
\label{def:covering-with-tail}
Let $\X$ be as in Def.~\ref{def:h-set-with-tail}. Let $N$, $M$
be h-sets with tails in $\X$ and $P : N \to \X$ be a continuous and compact mapping in $\X$.

We say that $N$ $P$-covers $M$
(denoted as before in Def.~\ref{def:singlevalued-covering-FINITE}
by $N \cover{P} M$), iff there exists a continuous and compact homotopy
$H : [0,1] \times |N| \to \X$ satisfying the conditions:
\begin{enumerate}
\item[C0.] $H\left(t, |N|\right) \subset \R^{d_M} \times |M_2|$;
\item[C1.] $H\left(0, \cdot\right) = P$;
\item[C2.] $H\left(\left[0,1\right], N_1^- \times |N_2|\right) \cap M = \varnothing$;
\item[C3.] $H\left(\left[0,1\right], |N|\right) \cap \left(M_1^+ \times |M_2|\right) = \varnothing$;
\item[C4.] there exists %a linear map $A : \R^{u} \to \R^{u}$ 
 $a\in\mathbb{R}\setminus[-1,1]$
and a point $\bar{r} \in M_2$ such that for all
$(p, q, r) \in N_c = [-1,1] \times \cB_{d_N-1}\times |N_2|$
we have:
\begin{eqnarray*}
H_c(1, (p, q, r)) & = & (a p, 0, \bar{r}), %\\
%A(\bd \Ball_u(0,1)) }& \subset & \R^u \setminus \overline{\Ball_u}(0,1)
\end{eqnarray*}
where again $H_c(t, \cdot) = c_M \circ H(t, \cdot) \circ c_N^{-1}$
is the homotopy expressed in good coordinates.
\end{enumerate}
\end{definition}
For some intuition of the infinite-dimensional covering relation, see Figure~\ref{fig:covrel}. 
The important fact is that such a  covering relation also fulfils the Itinerary Lemma, which was summarised in \cite{szczelina-zgliczynski-focm-2} in the following theorem:

\begin{theorem}[\cite{szczelina-zgliczynski-focm-2}]
\label{thm:periodic-covering-with-tail}
The claim of Theorem~\ref{thm:periodic-covering-FINITE} is true for
a covering relation loop \eqref{eq:loop} where sets $M_i$ are h-sets with tails in a
real Banach space $\X$, provided that the map $f$ is compact.
\end{theorem}

In practice, especially in computer-assisted proofs, we will use the following Lemma to check a covering relation with one unstable direction:

\begin{lemma}[\cite{szczelina-zgliczynski-focm-2}]
\label{lem:one-unstable}
For a h-set with tail $N$ define: 
\begin{itemize}
\item $N_c^{l} = \{-1\} \times \cB_{d_N-1} \times |N|$,\quad  $N^l = c_N^{-1}(N_c^{l})$ -- the \emph{left edge of $N$}, and 
\item $N_c^{r} = \{1\} \times \cB_{d_N-1} \times |N|$,\quad  $N^r = c_N^{-1}(N_c^{r})$ -- the \emph{right edge of $N$}.
\end{itemize}

Let $\X$ be a Banach space, $X \subset \X$ be an ANR,
$N = (N_1, |N_2|)$, $M = (M_1, |M_2|)$ be (horizontal) h-sets with tails in $X$
and $P : |N| \to X$ be a continuous and
compact map such that the following conditions apply
(with $P_c = c_M \circ P \circ c_N^{-1} : N_c \to M_c$):
\begin{enumerate}
\item[{\bf CC1.}]
	% in strip
	$\pi_{\X_2} P\left(|N|\right) \subset |M_2|$;
\item[{\bf CC2.}]
	% unstable boundary outside and unstable across
	Either {\bf CC2a}:\quad
$
\pi_1 P_c\left(N_c^l\right)\subset(-\infty,-1) 
\text{ \quad and \quad }
\pi_1 P_c\left(N_c^r\right) \subset(1,\infty) ,
$
\\
%	\begin{eqnarray*}
%	P_c\left(N_c^l\right) \subset (-\infty, -1) \times \R_{d_N-1} \times |M_2|
%	& and &
%	P_c\left(N_c^r\right) \subset (1, \infty) \times \R_{d_N-1} \times |M_2|
%	\end{eqnarray*}
	or {\bf CC2b}: \quad 
$
\pi_1 P_c\left(N_c^l\right)\subset(1,\infty) 
\text{ \quad and \quad }
\pi_1 P_c\left(N_c^r\right)\subset(-\infty,-1) .
$
%	\begin{eqnarray*}
%	P_c\left(N_c^l\right) \subset (1, \infty) \times \R_{d_N-1} \times |M_2|
%	& and &
%	P_c\left(N_c^r\right) \subset (-\infty, -1) \times \R_{d_N-1} \times |M_2|;
%	\end{eqnarray*}	
\item[{\bf CC3.}]
	% no touch entry set
	$P_c\left(N_c\right) \cap \left(M^+_c \times |M_2|\right) = \varnothing$
\end{enumerate}
Then $N \cover{P} M$ with the homotopy given as
$H(t, \cdot) = (1-t) \cdot P + t \cdot (A, 0, \bar{r})$,
where $A : \R \to \R$ such that $Ap = 2p$ (CC2a) or $Ap = -2p$ (CC2b)
and $\bar{r}$ is any selected point in $|M_2|$.
\end{lemma}

% SUBSECTION %%%%%%%%%%%%%%%%%%%%%%%%%%%%%%%%%%%%%%%%%%%%%%%%%%%%%%%%%%%%%%%%%%%
\subsection{Sharkovskii theorem via covering relations}

The goal of this work is to generalise this theorem, which is the main applicable result in the paper \cite{GZ2022}.
\begin{theorem}[Theorem~2 in \cite{GZ2022}]
\label{thm:sh-manyD}
Consider a set $G=I\times \cB(R)$, where $I\subset \mathbb{R}$ is a closed interval, and $\cB(R) \subset \mathbb{R}^{N-1}$ is a closed ($N-1$)-dimensional ball of radius $R$. Denote its points by $(p,q) \in I\times \cB(R)$.

Let $F:G \to \interior G$ be a continuous map with a periodic point $o_0=(p_0,q_0)$ with least period $n$, and denote its orbit by 
 \[
 \mathcal{O}=\{(p_i,q_i):=F^i(p_0,q_0), i=0,\dots,n-1\}.
 \]
Suppose that there exist $\delta_0$, $\delta_1$, \dots, $\delta_{n-1}>0$ such that the intervals $[p_i\pm\delta_i]\subset I$ are pairwise disjoint and
\[
\forall\: i\in\{0,\dots,n-1\} \qquad F\left([p_i\pm\delta_i]\times\cB(R)\right) \subset (p_{i+1}\pm\delta_{i+1})\times \B(R).
\]
Then for every natural number $m$ succeeding $n$ in the Sharkovskii order \eqref{eq:order} $F$ has a point with the least period $m$.
\end{theorem}

\begin{figure}[ht]
\begin{center}
	\includegraphics[width=0.9\linewidth]{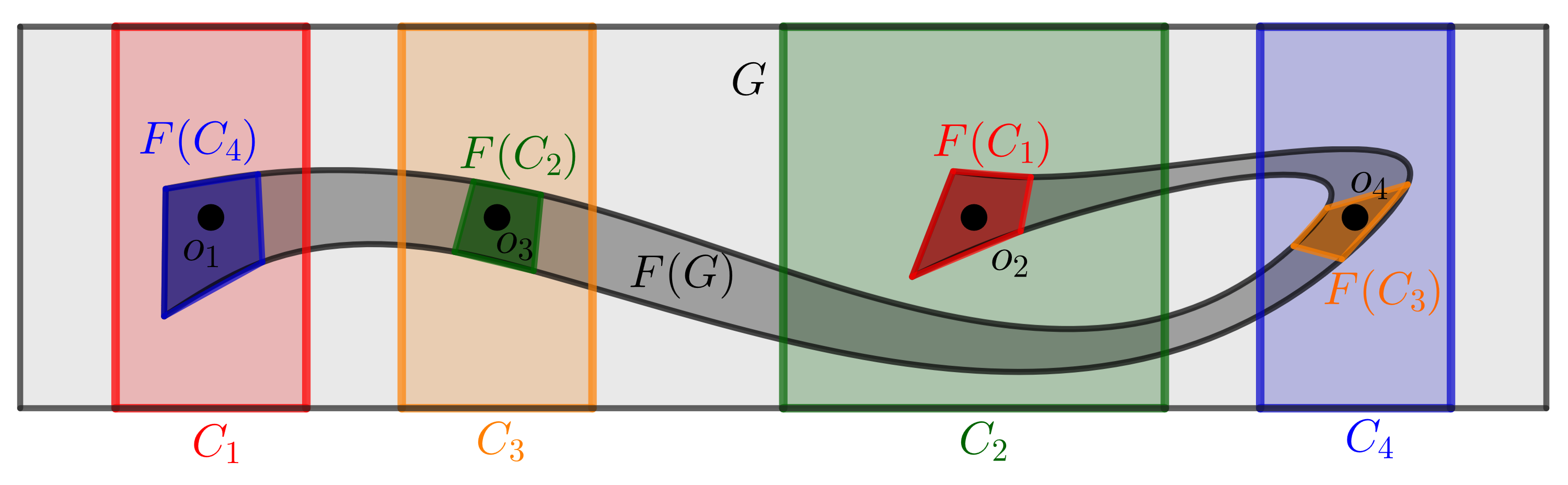}
	\caption{\label{fig:contraction}Illustration of assumptions of Theorem~\ref{thm:sh-manyD}:\newline $F(G)\subset\interior G$,\quad $F(C_{i})\subset\interior C_{i+1}$,\quad  $i=1,2,3,4$}
\end{center}
\end{figure}

For some intuition of the assumptions of Theorem \ref{thm:sh-manyD}, see Fig.~\ref{fig:contraction}.

\noindent\textit{Idea of the proof}: First the authors note that the map $F$ is `similar' to an one--dimensional map $f:J\to J$ on a closed interval $J$, which also has an $n$-periodic orbit $\mathcal{Q}$ with the same `permutation pattern'. For a fixed $m\triangleright n$,
they apply Theorem \ref{th:Burns} to the map $f$ and obtain a non-repeating $m$-loop of $\mathcal{Q}$-intervals $K_i$, proving the existence of an $m$-periodic point for $f$. Then, for every interval covering relation $K_i\overset{f}{\longrightarrow}K_{i+1}$ from this loop, they construct h-sets $S(K_i)$, $S(K_{i+1})$, such that $S(K_i)\overset{F}{\Longrightarrow} S(K_{i+1})$ and therefore they obtain a corresponding non-repeating loop of horizontal covering relations:
\begin{equation*}
	  S(K_0)
	  \overset{F}{\Longrightarrow}
	  S(K_1)
	  \overset{F}{\Longrightarrow}
	  \dots
	  \overset{F}{\Longrightarrow}
	  S(K_{m-1})
	  \overset{F}{\Longrightarrow}
	  S(K_0)\text{,}
	\end{equation*}
which proves the existence of a periodic orbit with fundamental period $m$ also for $F$.\qed

The above theorem is easily applicable to finite dimensional maps with an attracting periodic orbit and one unstable direction, also with computer assistance. If we want to prove a Sharkovskii-type theorem for a map $F$ with an attracting $n$-periodic orbit $\O = \{o_0,\dots,o_i=F^i(o_0),\dots,o_n=o_0\}$, it is sufficient to find $n+1$ cube-like sets:
\begin{itemize}
    \item a neighbourhood of the whole orbit $G\supset \O$ such that $F(G)\subset \interior G$, and
    \item neighbourhoods $C_i\ni o_i$ such that $F(C_i)\subset \interior C_{i+1}$ for every orbit point $o_i$.
\end{itemize}

\begin{figure}[ht]
	\begin{center}\includegraphics[width=0.605\linewidth]{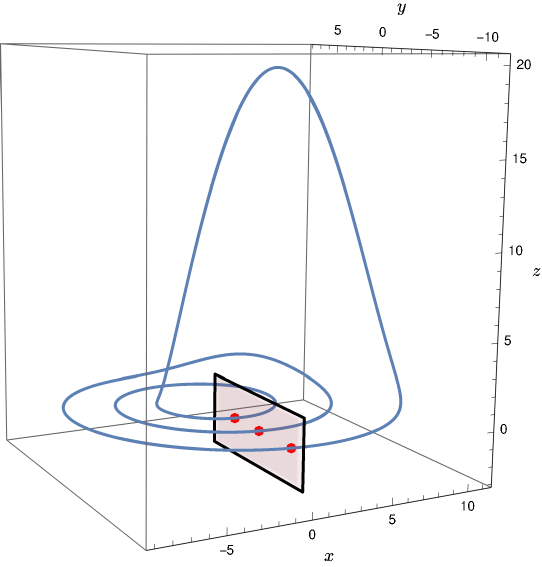}\end{center}
	\caption{\label{fig:attr3}The attracting $3$-periodic orbit for the system \eqref{eq:rossler525}}
\end{figure}

As an application, recall the toy example from \cite{GZ2022}, the R\"ossler system with $a=5.25$ and $b=0.2$:
\begin{equation}\label{eq:rossler525}
\begin{cases}
x'=-y-z,
\\
y'=0.2 y+x,
\\
z'=z (x-5.25)+0.2
\end{cases}\text{,}
\end{equation}
which was proven to have an attracting 3-periodic orbit \cite{GZ2021}, see Fig.\ \ref{fig:attr3}. 

This example will be used later, as a base for the proof
of similar behaviour in infinite-dimensional system, so we would like
to re-state elements of the computer assisted proof from \cite{GZ2022}.

\begin{figure}[ht]
\begin{center}
	\includegraphics[width=0.8\textwidth]{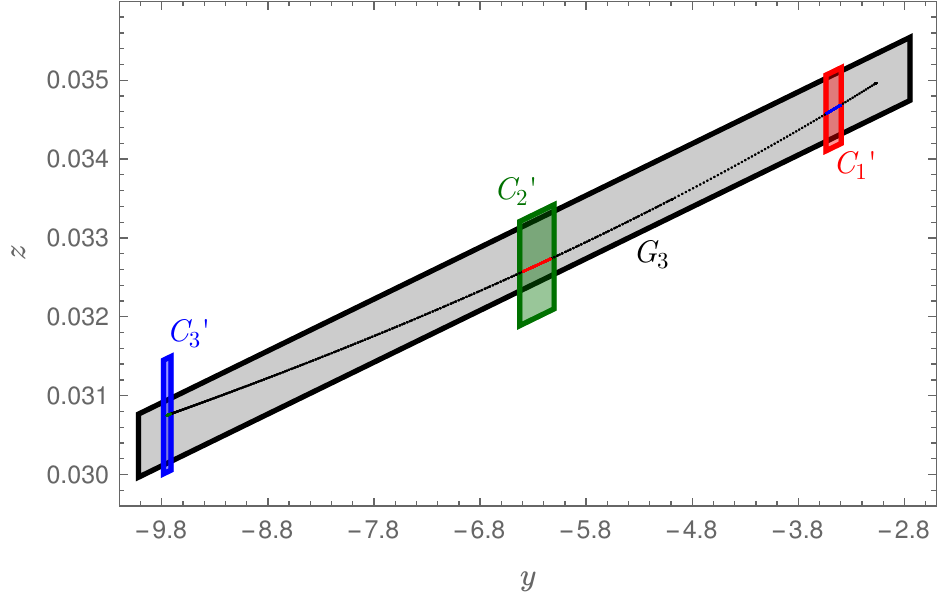}
	\caption{\label{fig:grid3}The orbit neighbourhood $G_3$ (grey) from Lemma \ref{lem:Roessler525} and its image through $P$ (black, thin `curve' inside $G_3$).
	\newline The supersets $C_i'$ of orbit points' neighbourhoods, $i=1,2,3$, and their images are marked in red, green and blue.}
\end{center}
\end{figure}

The orbit is $3$-periodic for the Poincar\'e map $P$ of the system, partly defined on the section $\Pi=\{x=0,y<0\}$ and the map $P$ exhibits a strong contraction in the $z$ direction, which means that it can be seen as a perturbation of  a one-dimensional map and Theorem \ref{thm:sh-manyD} should apply.
Indeed, the following lemma was proven by checking the inclusions with the use of \texttt{C++} library \texttt{CAPD} \cite{capd-article}:

\begin{lemma}[\cite{GZ2022}]\label{lem:Roessler525}
Let $M=\left[\smallmatrix -1.& 0.000656767 \\ -0.000656767 & -1. \endsmallmatrix\right]$. The parallelogram $G_3$ in the $(y, z)$ coordinates on the Poincar\'e section $\Pi$ (see
Fig. \ref{fig:grid3}):
\[
G_3= \bmatrix -6.38401 \\ 0.0327544 \endbmatrix
+ M \cdot
\bmatrix \pm 3.63687 \\ \pm 0.0004 \endbmatrix
%\text{,}
\]
fulfils the assumptions of Theorem \ref{thm:sh-manyD} for the Poincar\'e map $P$ with the attracting $3$-periodic orbit from \cite[Lemma 4]{GZ2021}. The orbit points neighbourhoods $C_i$, $i=1,2,3$ are defined by $C_i = C'_i \cap G_3$, where
\begin{align*}
	&C'_1= \bmatrix -3.46642 \\ 0.0346316 \endbmatrix
+ M \cdot
\bmatrix \pm 0.072 \\ \pm 0.00048 \endbmatrix
	\text{,}
	\\
	&C'_2= \bmatrix -6.26401 \\ 0.0326544 \endbmatrix
+ M \cdot
\bmatrix \pm 0.162 \\ \pm 0.00066 \endbmatrix
	\text{,}
	\\
	&C'_3= \bmatrix -9.74889\\ 0.0307529 \endbmatrix
+ M \cdot
\bmatrix \pm 0.036 \\ \pm 0.00072 \endbmatrix.
\end{align*}
Therefore, for every $n\in\N$, the system \eqref{eq:rossler525} has an $n$-periodic orbit for $P$, passing through $G$.
\end{lemma}

Let us now formulate the generalisation of Theorem~\ref{thm:sh-manyD} for a 
compact map over possibly infinite-dimensional space, which has an attracting $n$-periodic orbit. We require, as in the case of Theorem~\ref{thm:sh-manyD}, 
that essential dynamic of the map is in some sense close to 1-dimensional map. 
\begin{theorem}
\label{thm:sh-infD}
Let $\mathcal{X} = \mathcal{X}_1 \times \mathcal{X}_2$ be
some Banach spaces as in Definition~\ref{def:h-set-with-tail},
with $\mathcal{X}_1 = \R^N$.

Consider a set $G=I\times \cB(R) \times \Xi$ $\subset \mathcal{X}$, 
where $I\subset \mathbb{R}$ is a closed interval, 
$\cB(R) \subset \mathbb{R}^{N-1}$ is a closed ball of radius 
$R$ and $\Xi$ is a closed, convex and bounded subset of the 
Banach space $\X_2$.

Let $F:G \to \interior I\times \B(R)\times \Xi$ be a continuous and compact map with a periodic point $o_0=(p_0,q_0,r_0) %\in \interior I\times \cB(R)\times \Xi
$ with least period $n$, and denote its orbit by 
\[
\mathcal{O}=\{o_i=(p_i,q_i,r_i):=F^i(p_0,q_0,r_0),\: i=0,\dots, n-1\}%\subset \interior I\times \B(R)\times \Xi
\text{.}
\]
Suppose also  that there exist $\delta_0$, $\delta_1$, \dots, $\delta_{n-1}>0$ such that the intervals $[p_i\pm\delta_i]\subset I$ are pairwise disjoint and
\begin{equation}
\label{eq:contraction}
\forall\: i\in\{0,\dots,n-1\} \quad F\left([p_i\pm\delta_i]\times\cB(R)\times \Xi\right) \subset (p_{i+1}\pm\delta_{i+1})\times \B(R)\times \Xi.
\end{equation}
Then for every natural number $m$ succeeding $n$ in the Sharkovskii order \eqref{eq:order} $F$ has a point with the least period $m$.
\end{theorem}

\noindent\textit{Proof}: 
We shall follow the proof of Theorem \ref{thm:sh-manyD}, or \cite[Thm.\ 2]{GZ2022}, but we want to avoid the \emph{contracting grid} formality introduced there. 

First, note that the assumptions set the periodic orbit $\O$ in some order from the point of view of the first coordinate $p\in I$. Therefore, it is useful to 
construct a `model' continuous map $f: J\to J$ on some interval $J$ which also has an $n$-periodic orbit with the same permutation pattern. 
In what follows we will use projection $\pi_1$ to denote
the projection onto the first coordinate of the finite-dimensional 
space $\mathcal{X}_1 = \R^N=\R\times\R^{N-1}$,
that is, the projection onto the set $I$.

Such a `model' map can be defined, for instance, piecewise linearly. Let $J=[\min\pi_1(\O),\max\pi_1(\O)]\subset I$ and then $f:J\to J$ is defined as follows
\[
f(p)= 
\begin{cases}
p_{i+1} &\text {for } p=p_i,
\\[2mm]
p_{j+1}+\frac{p-p_j}{p_k-p_j}(p_{k+1}-p_{j+1}) &\text{for $p\in (p_j,p_k)$ such that}
\\[-2mm]
&\pi_1(\O)\cap(p_j,p_k)=\varnothing.
\end{cases}
\]
Then, the set $\mathcal{Q}=\pi_1(\O) = \{p_0,\dots,p_{n-1}\}$ is the $n$--periodic orbit for $f$ with the same pattern as $\O$, in the sense that, for all $i=0,\dots,n-1$:
\[
f(\pi_1(o_i))=f(p_i) = p_{i+1} = \pi_1(o_{i+1})=\pi_1(F(o_i)).
\]

Now, let us fix $m\triangleright n$. From Theorem \ref{th:Burns} we know that there exists a non-repeating %$\O$-forced 
	loop of $m$ $\mathcal{Q}$-intervals $J_i$, $i=0,\dots, m-1$:
	\begin{equation}\label{eq:interval_loop}
	  J_0  \overset{f}{\longrightarrow} J_1 \overset{f}{\longrightarrow} J_2 \overset{f}{\longrightarrow} \dots \overset{f}{\longrightarrow} J_{m-1} \overset{f}{\longrightarrow} J_0\text{,}
	\end{equation}
	which proves the existence of an $m$-periodic point for $f$ in the interval $J$.

\begin{figure}[ht]
	\includegraphics[width=0.99\linewidth]{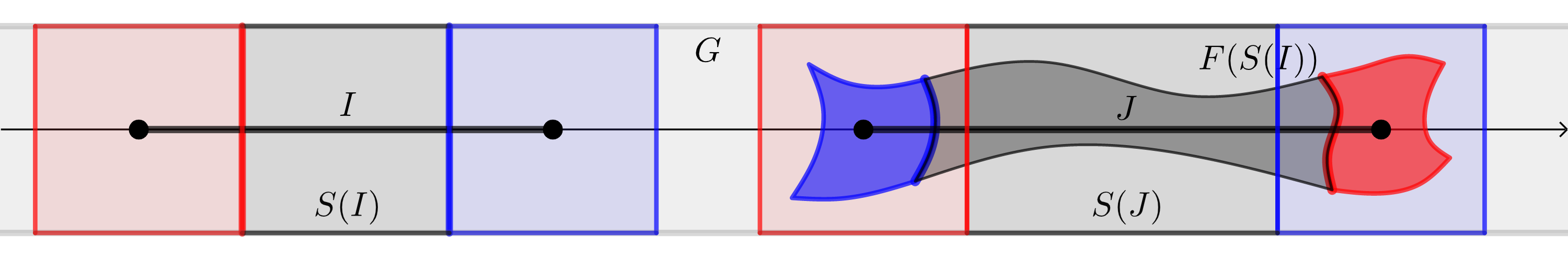}
	\caption{\label{fig:S(J)}Construction of $S(I)$, $S(J)$ (medium grey) over the intervals $I$, $J$. The neighbourhoods of the interval ends (red and blue) are mapped to the interiors of consecutive neighbourhoods which forces the image $F(S(I))$ to be spanned across $S(J)$ and, as a consequence, the covering relation $S(I) \overset{F}{\Longrightarrow} S(J)$.}
\end{figure}

Next, we want to construct a similar loop of horizontal covering relations for $F$. To obtain this, first define the corresponding h-sets with tails. Consider any interval $J_i$ from the loop \eqref{eq:interval_loop} and suppose that its ends are the two (not necessarily consecutive) points $p_a<p_b$ from the periodic orbit. Then, the corresponding h-set with tail $S(J_i)$ will be defined as
\[
\begin{cases}
|S(J_i)|&=\left([p_a+\delta_a,p_b-\delta_b]\times \cB(R)\right)\times \Xi\text{,}
\\[2mm]
c_{S(J_i)}(p,q,r)&=\left(-1+2\dfrac{p-p_a-\delta_a}{p_b-\delta_b-p_a-\delta_a},q,r\right).
\end{cases} 
\]
For some intuition of construction of the head of $S(J_i)$, see Fig.\ \ref{fig:S(J)}. Please note that $\O \cap |S(J_i)| = \varnothing$. 

We claim that
\begin{lemma}
\label{lem:h-sets-tails-covering}
The h-sets with tails $S(J_i)$, $j=0,\dots,m-1$ fulfil the loop of horizontal covering relations
\begin{equation}\label{eq:Final_Loop}
	  S(J_0)
	  \overset{F}{\Longrightarrow}
	  S(J_1)
	  \overset{F}{\Longrightarrow}
	  S(J_2)
	  \overset{F}{\Longrightarrow}
	  \dots
	  \overset{F}{\Longrightarrow}
	  S(J_{m-1})
	  \overset{F}{\Longrightarrow}
	  S(J_0)\text{.}
	\end{equation}
\end{lemma}

\noindent\textit{Proof}: Fix any $i=0,\dots,m-1$ and consider the interval covering relation $J_i \overset{f}{\longrightarrow} J_{i+1}$ from the loop \eqref{eq:interval_loop}. Denote the ends of $J_i$ by $p_a<p_b$ and those of $J_{i+1}$ by $p_c<p_d$. From the definition we know that
\[
\min \{f(p_a), f(p_b)\} \leq p_c \text{ \qquad and \qquad } p_d \leq \max \{f(p_a),f(p_b)\}.
\]
Let us check first the case $f(p_a)\leq p_c < p_d \leq f(p_b)$, or, in other words $p_{a+1}\leq p_c < p_d \leq p_{b+1}$ Then we know that for $F$ also
\[
\pi_1 (o_{a+1})=\pi_1 (F(o_a))\leq \pi_1(o_c) < \pi_1(o_d) \leq \pi_1 (F(o_b))=\pi_1 (o_{b+1}).
\]

Now, for the sets $S(J_i)$ and $S(J_{i+1})$ we check the three conditions from Lemma \ref{lem:one-unstable}:

\begin{enumerate}
\item[{\bf CC1.}]	% in strip
	$\pi_{\X_2} F\left(|S(J_i)|\right) \subset \Xi$, because  $F:G \to \interior \left(I\times \cB(R)\right)\times \Xi$.
\item[{\bf CC2a.}] % unstable boundary outside and unstable across
The left edge:
    \begin{multline*}
        \pi_1 F_c\left(S(J_i)_c^l\right)=\pi_1 c_{S(J_{i+1})}(F(S(J_i)^l))=
        \\
        =\pi_1 c_{S(J_{i+1})}(F(\{p_a+\delta_a\}\times\cB(R)\times \Xi))\overset{\eqref{eq:contraction}}{\subset}
        \\
        \subset\pi_1 c_{S(J_{i+1})}((p_{a+1}\pm\delta_{a+1})\times\cB(R)\times \Xi)=
        \\
        =\left(-1+2\frac{p_{a+1}-\delta_{a+1}-p_c-\delta_c}{p_d-\delta_d-p_c-\delta_c},
        -1+2\frac{p_{a+1}+\delta_{a+1}-p_c-\delta_c}{p_d-\delta_d-p_c-\delta_c}\right)\subset 
        \\
        \subset (-\infty,-1)\text{,}
    \end{multline*}
because either $p_{a+1}<p_c$ and then $p_{a+1}+\delta_{a+1}<p_c+\delta_c$ (as the $[p_i\pm\delta_i]$ intervals are disjoint by assumption),
or $p_{a+1}=p_c$ and then $\delta_{a+1}=\delta_c$, so, in general, $p_{a+1}+\delta_{a+1}\leq p_c+\delta_c$.

The right edge similarly:
    \begin{multline*}
        \pi_1 F_c\left(S(J_i)_c^r\right)=\pi_1 c_{S(J_{i+1})}(F(S(J_i)^r))=
        \\
        =\pi_1 c_{S(J_{i+1})}\left(F\left(\{p_b+\delta_b\}\times\cB(R)\times \Xi\right)\right)\overset{\eqref{eq:contraction}}{\subset}
        \\
        \subset\pi_1 c_{S(J_{i+1})}\left((p_{b+1}\pm\delta_{b+1})\times\cB(R)\times \Xi\right)=
        \\
        =\left(-1+2\frac{p_{b+1}-\delta_{b+1}-p_c-\delta_c}{p_d-\delta_d-p_c-\delta_c},
        -1+2\frac{p_{b+1}+\delta_{b+1}-p_c-\delta_c}{p_d-\delta_d-p_c-\delta_c}\right)\subset 
        \\
        \subset (1,\infty)\text{,}
    \end{multline*}
    which follows from $p_d-\delta_d \leq p_{b+1}-\delta_{b+1}$.
\item[{\bf CC3.}]
	% no touch entry set
	$F_c\left(S(J_i)_c\right) \cap \left(S(J_{i+1})^+_c \times \Xi\right) = \varnothing$, because $F:G \to \interior \left(I\times \cB(R)\right)\times \Xi$.
\end{enumerate}
The case $f(p_b)\leq p_c < p_d \leq f(p_a)$ follows analogously, but we check the condition \textbf{CC2b}. This finishes the proof of Lemma~\ref{lem:h-sets-tails-covering}.
\qed

Now, to finish the proof of Theorem~\ref{thm:sh-infD}, note that since the loop \eqref{eq:interval_loop} was non-repeating, then so is the corresponding loop \eqref{eq:Final_Loop}, and we deduce that 
\[
%\interior 
S(J_0) \cap \bigcup_{i=1}^{m-1} S(J_{i}) = \varnothing\text{,}
\]
which means that the $m$-periodic orbit for $F$ which existence follows from the loop \eqref{eq:Final_Loop} has fundamental period $m$.\qed

\begin{rem}
Note that in the proof of Theorem \ref{thm:sh-infD} we do not need the assumption of existence of the $n$-periodic orbit $\O$. Its existence, in fact, follows from the relations between the sets $[p_i\pm\delta_i]\times \cB(R)\times \Xi$ and their images, and Schauder Fixed Point Theorem.
\end{rem}

\begin{definition}\label{def:SharProp}
If some restriction of a map $f:\mathcal{X}\to\mathcal{X}$ on a Banach space fulfils the assumptions of the above Theorem \ref{thm:sh-infD} (or, in the finite-dimensional case, Thm.\ \ref{thm:sh-manyD}), we will say that $f$ has \emph{Sharkovskii property}. Less formally, we will also say that a flow has Sharkovskii property if it admits a Poincar\'e map with this property.
\end{definition}

%% file: 3_perturbation.tex
%%%%%%%%%%%%%%%%%%%%%%%%%%%%%%%%%%%%%%%%%%%%%%%%%%%%%%%%%%%%%%%%%%%%%%%%%%%%%%%%
% SECTION %%%%%%%%%%%%%%%%%%%%%%%%%%%%%%%%%%%%%%%%%%%%%%%%%%%%%%%%%%%%%%%%%%%%%%
%%%%%%%%%%%%%%%%%%%%%%%%%%%%%%%%%%%%%%%%%%%%%%%%%%%%%%%%%%%%%%%%%%%%%%%%%%%%%%%%
\section{\label{sec:perturbation}Delayed perturbation of a chaotic ODE}

In this section we apply Theorem~\ref{thm:sh-infD} to a general
scenario of an ODE perturbed by a small term in a form of Delay Differential 
Equation (DDE). More formally, we consider an autonomous ODE in $\R^d$:
\begin{equation}
\label{eq:ode}
u'(t) = f\left(u(t)\right)
\end{equation}
with initial data $u(0) = u_0$ and its perturbation by a small term:
\begin{equation}
\label{eq:dde}
u'(t) = f\left(u(t)\right) + \epsi \cdot g\left(u(t - \tau)\right)
\end{equation}
where the delay $\tau$ is fixed and reasonably big, e.g. $\tau = 1$ 
(contrary to the work \cite{wojcik-zgliczynski-dde-chaos} where the delay was 
sufficiently small, $\tau \approx \epsi$). The initial data is continuous
$u(s) = \phi(s)$ for $s \in [-\tau, 0]$, as is commonly used in DDEs. Since we 
are interested in the dynamical behaviour of \eqref{eq:dde} in comparison 
to \eqref{eq:ode}, we usually set $\phi(0) = u_0$, but we will always state 
it explicitly. We will assume that the ODE~\eqref{eq:ode} has the Sharkovskii 
property as in Theorem~\ref{thm:sh-manyD}, and then we will show that a 
sufficiently small (in $\epsi$) perturbation preserves this property. 
The delay here might be arbitrary large (but constant). The small $\epsi$ 
will depend on $f, g$ and $\tau$. 

% SUBSECTION %%%%%%%%%%%%%%%%%%%%%%%%%%%%%%%%%%%%%%%%%%%%%%%%%%%%%%%%%%%%%%%%%%%
\subsection{\label{sec:dde-glossary}Assumptions and glossary for DDEs}

From now on, we assume that a fixed $\tau > 0$ is given. 

\begin{itemize}
\item 
Recall that by $\| \cdot \|$ we denote the maximum norm on $\R^d$.
We will denote by $\mathcal{C} = \mathcal{C}\left([-\tau, 0], \R^d\right)$
the set of all continuous functions from interval $[-\tau, 0]$ to $\R^d$.
We endow $\mathcal{C}$ with the standard supremum norm 
derived from the chosen norm in $\R^d$ and denote it by:
\begin{equation*} 
\| \phi \|_\infty = \max_{t \in [-\tau, 0]} \| \phi(t) \|
\end{equation*}
so that $\mathcal{C}$ is a Banach space.

\item
We will use the convention that by $u$ we denote the solution of \eqref{eq:ode} 
and by $\bar{u}$ the solution of \eqref{eq:dde}. 

\item
Let $u : \R \to \R^d$ be any function. By $u_t \in \mathcal{C}$ we will 
denote the \emph{segment} of $u$ at $t$, i.e. $u_t(s) = u(s + t)$ 
for $s \in [-\tau, 0]$.

\item
Let $K \subset \R^d$. Following the notation from \cite{wojcik-zgliczynski-dde-chaos},
by $K_{\mathcal{C}}$ we denote the set $\{ \phi \in \mathcal{C} : \phi(0) \in K \}$ 
and by $K_{\mathcal{C}}(\epsi)$ we denote the set 
$\{ \phi \in K_\mathcal{C} : \| \phi - \phi(0) \|_\infty \le \epsi \}$.
Next, for a given constant $L$ we denote by $K_{\mathcal{L}(L)}$ the set of 
Lipschitz functions $\phi$ with constant $L$
with $\phi(0) \in K$, i.e.:
\begin{equation*}
K_{\mathcal{L}(L)} = \left\{ \phi \in K_\mathcal{C} : \forall_{t_1, t_2 \in [-\tau, 0]} \| \phi(t_2) - \phi(t_1) \| \le L |t_2 - t_1 | \right\}.
\end{equation*}
Denote $\Omega = \{ 0 \} \subset \R^d$, then 
$\Omega_{\mathcal{L}(L)}$ is a set of Lipschitz functions $\phi \in \mathcal{C}$
with Lipschitz constant $L$ and the value $\phi(0) = 0$.
In $K_{\mathcal{L}(L)}$ we introduce the following coordinates for 
$\phi$: $(\phi(0), \eta)$, where $\eta(s) = \phi(s) - \phi(0)$ for 
$s \in [-\tau, 0]$. Now, we can write 
$K_{\mathcal{L}(L)} = K \times \Omega_{\mathcal{L}(L)}$.

\item
Thorough the article we assume the following condition:

\textbf{Condition A}\label{condA}: $f$ is continuous, bounded by 
$M_f > 0$ ($f(u) \le M_f$ for all $u \in \R^d$)
and Lipschitz with the constant 
$L_f > 0$ ($\|f(u_1) - f(u_2)\| \le L_f \|u_1 - u_2\|$), and
the same for $g$, with constants $M_g$ and $L_g$, respectively.
\emph{Note that from now on, we will not state those assumptions 
in Theorems/Lemmas explicitly, but we might remind about it from time to time.}

\item
We will follow the convention that $K_{\mathcal{L}}$ denotes 
the set $K_{\mathcal{L}(M_f + \epsi M_g)}$.
\end{itemize}

\subsection{Basic estimates}

First, restate the 
equations \eqref{eq:ode} and \eqref{eq:dde} as the integral
equations:
\begin{equation}
\label{eq:ode-integral}
u(t) = u(0) + \int_0^t f(u(s)) ds
\end{equation}
\begin{equation}
\label{eq:dde-integral}
\bar{u}(t) = \bar{u}(0) + \int_0^t f(\bar{u}(s)) ds + \epsi \int_0^t g(\bar{u}(s - \tau)) ds
\end{equation}

Basic estimates for solutions $u$, $\bar{u}$ will follow from the 
classic Gr\"onwall integral inequality:
\begin{lemma}[Gr\"onwall inequality]
Let $I = [a, b]$ or $I = [a,\infty)$ or $I = [a, b)$.
Let $u : I \to \R$ be such that there exist continuous functions
$\alpha, \beta : I \to \R$,
with $\beta(t) \ge 0$ for all $t \in I$, $\alpha$ non-decreasing,
such that
\begin{equation}
u(t) \le \alpha(t) + \int_a^t \beta(s) u(s) ds.
\end{equation}
Then
\begin{equation}
u(t) \le \alpha(t) \exp \left(\int_a^t \beta(s) ds\right),
\end{equation}
for all $t \in I$.
\end{lemma}

We start with a couple of estimates which will be used later. 
Obviously we have:
\begin{lemma}
\label{lem:t-lipshitz}
If $\bar{u}$ is a solution to \eqref{eq:dde} then 
for any $t_1, t_2 \ge 0$ we have 
\begin{equation}
\| \bar{u}(t_2) - \bar{u}(t_1)\| \le (M_f + |\epsi| M_g) \cdot |t_2 - t_1|.
\end{equation}
\end{lemma}
\noindent\textit{Proof}: W.l.o.g. assume $t_2 \ge t_1$. 
Using \eqref{eq:dde-integral} we get:
\begin{eqnarray*}
\| \bar{u}(t_2) - \bar{u}(t_1) \| 
	& = & \left\| \int_{t_1}^{t_2} f(\bar{u}(s)) ds + \epsi \int_{t_1}^{t_2} g(\bar{u}(s - \tau)) ds \right\|  \\
	& \le & (M_f + |\epsi| M_g) \cdot |t_2 - t_1|
\end{eqnarray*}
\qed

\begin{lemma}
\label{lem:ode-dde-delta}
If $u$ is a solution of \eqref{eq:ode} and $\bar{u}$ of \eqref{eq:dde} 
with $u(0) = \bar{u}(0)$ then:
\begin{equation}
\label{eq:ode-dde-delta}
\| u(t) - \bar{u}(t) \| \le |\epsi| \cdot M_g \cdot t \cdot e^{L_f t}
\end{equation}
\end{lemma}
\noindent\textit{Proof}: using \eqref{eq:ode-integral} 
and \eqref{eq:dde-integral} we get:
\begin{eqnarray*}
\| u(t) - \bar{u}(t) \| 
	& = & \left\| \int_0^t f(u(s)) - f(\bar{u}(s)) ds + \epsi \int_0^t g(\bar{u}(s - \tau)) ds \right\|  \\
	& \le & \int_0^t \|f(u(s)) - f(\bar{u}(s))\| ds + |\epsi| \int_0^t \|g(\bar{u}(s - \tau))\| ds \\
	& \le & L_f \int_0^t \|u(s)) - \bar{u}(s)\| ds + |\epsi| M_g t
\end{eqnarray*}
Now \eqref{eq:ode-dde-delta} follows from Gr\"onwall's inequality.
\qed

\begin{lemma}
\label{lem:dde-dde-delta}
If $\bar{v}, \bar{u}$ be solutions of \eqref{eq:dde}, then
\begin{eqnarray*}
\| \bar{v}(t) - \bar{u}(t) \| 
	& \le & |\epsi| \cdot L_g \cdot \tau \cdot e^{\left(L_f + |\epsi| L_g\right) \cdot t} \cdot \| \bar{v}_0 - \bar{u}_0 \|_\infty  \ + \ \\
	& + & \| \bar{v}(0) - \bar{u}(0) \| \cdot e^{\left(L_f + |\epsi| L_g\right) \cdot t}
\end{eqnarray*}
\end{lemma}
\noindent\textit{Proof}: let denote 
$w(t) = \| \bar{v}(t) - \bar{u}(t) \|$, then we have:
\begin{eqnarray*}
w(t) & = & \| \bar{v}(t) - \bar{u}(t) \| \\
	& \le & \| \bar{v}(0) - \bar{u}(0) \| + \int_0^t \| f(\bar{v}(s)) - f(\bar{u}(s)) \| ds + \epsi \int_0^t \| g(\bar{v}(s - \tau)) - g(\bar{u}(s - \tau)) \| ds  \\
	& \le & \| \bar{v}(0) - \bar{u}(0) \| + L_f \int_0^t w(s) ds + |\epsi| L_g \int_0^t w(s - \tau) ds \\
	& \le & \| \bar{v}(0) - \bar{u}(0) \| + L_f \int_0^t w(s) ds + |\epsi| L_g \int_{-\tau}^{t-\tau} w(s) ds \\
	& \le & \| \bar{v}(0) - \bar{u}(0) \| + L_f \int_0^t w(s) ds + |\epsi| L_g \int_{0}^{t} w(s) ds + |\epsi| L_g \int_{-\tau}^{0} w(s) ds \\
	& \le & \| \bar{v}(0) - \bar{u}(0) \| + \left(L_f + |\epsi| L_g\right) \int_0^t w(s) ds + |\epsi| L_g \tau \max_{s \in [-\tau, 0]} w(s) ds \\
	& = & \left(L_f + |\epsi| L_g\right) \int_0^t w(s) ds + |\epsi| L_g \tau \| \bar{v}_0 - \bar{u}_0 \|_\infty + \| \bar{v}(0) - \bar{u}(0) \|\text{,}
\end{eqnarray*}
and applying again Gr\"onwall's Lemma, we get the thesis.
\qed

From the proof of Lemma~\ref{lem:dde-dde-delta} we can prove even more:
\begin{cor}
If $\bar{v}, \bar{u}$ are solutions of \eqref{eq:dde} then:
\begin{itemize}
\item 
$\| \bar{v}(t) - \bar{u}(t) \| \le \left(1 + |\epsi| \cdot L_g \cdot \tau \right) \cdot e^{\left(L_f + |\epsi| L_g\right) \cdot t} \cdot \| \bar{v}_0 - \bar{u}_0 \|_\infty$

\item
if $\| \bar{v}(0) - \bar{u}(0)\| \le D \cdot |\epsi| \cdot L_g \cdot \tau \| \bar{v}_0 - \bar{u}_0 \|_\infty$ then 
\begin{equation*}
\| \bar{v}(t) - \bar{u}(t) \| \le |\epsi| \cdot (D + 1) \cdot  L_g \cdot \tau \cdot e^{\left(L_f + |\epsi| L_g\right) \cdot t} \cdot \| \bar{v}_0 - \bar{u}_0 \|_\infty
\end{equation*}

\item
in particular, if $\bar{v}(0) = \bar{u}(0)$ then 
\begin{equation*}
\| \bar{v}(t) - \bar{u}(t) \| \le |\epsi| \cdot  L_g \cdot \tau \cdot e^{\left(L_f + |\epsi| L_g\right) \cdot t} \cdot \| \bar{v}_0 - \bar{u}_0 \|_\infty
\end{equation*}
and, in turn:
\begin{equation}
\label{eq:small-tail}
\| \bar{v}_t - \bar{u}_t \|_\infty \le |\epsi| \cdot  L_g \cdot \tau \cdot e^{\left(L_f + |\epsi| L_g\right) \cdot t} \cdot \| \bar{v}_0 - \bar{u}_0 \|_\infty.
\end{equation}
\end{itemize}
\end{cor}

Another variant of Lemma~\ref{lem:dde-dde-delta}, that
does not need assumption on $g$ being Lipschitz:
\begin{lemma}
\label{lem:dde-dde-delta-2}
If $\bar{v}, \bar{u}$ be solutions of \eqref{eq:dde}, then
\begin{eqnarray*}
\| \bar{v}(t) - \bar{u}(t) \| 
	& \le & 2 \cdot |\epsi| \cdot M_g \cdot t \cdot e^{L_f t} \\ 
	& + & e^{L_f \cdot t} \cdot \| v(0) - u(0) \|.
\end{eqnarray*}
\end{lemma}
\noindent\textit{Proof}: 
let $v, u$ be the solutions of the ODE \eqref{eq:ode} with $u(0) = \bar{u}(0)$
and $v(0) = \bar{v}(0)$. Then we have:
\begin{eqnarray*}
\| \bar{v}(t) - \bar{u}(t) \|
	& = & \| \bar{v}(t) - v(t) + v(t) - u(t) + u(t) - \bar{u}(t) \|  \\
	& \le & \| \bar{v}(t) - v(t) \| + \| v(t) - u(t) \| + \| u(t) - \bar{u}(t) \|  \\
	& \le & 2 \cdot |\epsi| \cdot M_g \cdot t \cdot e^{L_f t} + e^{L_f t} \cdot \| v(0) - v(0) \|.
\end{eqnarray*}
Here we used Lemma~\ref{lem:ode-dde-delta} which also does not use Lipschitz 
assumption on $g$, and a standard estimate for ODEs for the middle term, e.g.
apply Lemma~\ref{lem:dde-dde-delta} with $\epsi = 0$).
\qed

From the proof of Lemma~\ref{lem:dde-dde-delta-2} we also have
a convenient formula:
\begin{pro}
Let $v, u$ be solutions of ODE~\eqref{eq:ode} and $\bar{v}, \bar{u}$
corresponding solutions of a DDE~\eqref{eq:dde} with $v(0) = \bar{v}(0)$,
$u(0) = \bar{u}(0)$. Then:
\begin{eqnarray*}
\| \bar{v}(t) - \bar{u}(t) \|
	& \le & 2 \cdot |\epsi| \cdot M_g \cdot t \cdot e^{L_f t} + \| v(t) - u(t) \|.
\end{eqnarray*}
\end{pro}

In all the formulas we are happy when $\epsi$ is factorised from an expression, 
because we can make it as small as we want in this case. More formally, for any $\delta > 0$
we can select $\epsi_M > 0$ such that for all $|\epsi| \le \epsi_M$ the solution $\bar{u}$
of the DDE~\eqref{eq:dde} will be $\delta$-close to the solution $u$ of the 
ODE~\eqref{eq:ode} with explicitly given estimates. 
We will rework all the estimates in the following section.

We will use those estimates in a similar way it was done in \cite{wojcik-zgliczynski-dde-chaos}. 
We will do one small change. Instead of using sets $K_\mathcal{C}$ 
from the mentioned work, we will
use $K_{\mathcal{L}} \subset K_{\mathcal{C}}$ as defined
in Section~\ref{sec:dde-glossary}. 
This will allow for some improvements:
\begin{enumerate}
\item 
	if $K$ is compact then $K_{\mathcal{L}}$ is compact 
	in $\mathcal{C}$ (essentially by Arzela-Ascoli theorem)

\item 
	convex and compact $K_{\mathcal{L}} = K \times \Omega_\mathcal{L}$ 
    is an Absolute Neighbourhood Retract (ANR), see
	\cite[\S~11 Corollary (4.4)]{granas-book}. It is important because the fixed point
	index is well defined for ANRs and has the usual properties.  

\item 
	if $\bar{u}$ is a solution to DDE with 
	$\bar{u}_0 \in K \times \Omega_\mathcal{L}$ then 
	$\bar{u}_t \in \tilde{K} \times \Omega_\mathcal{L}$ 
	for some other compact $\tilde{K}$
	\emph{for any $t$}. This is in contrast to 
	\cite{szczelina-zgliczynski-focm-1, szczelina-zgliczynski-focm-2} 
	where we had to assume $t \ge \tau$ to define proper Poincar\'e maps.
	Moreover, we do not need to assume $t \ge \tau$ to get a compact mapping. 
	The sets are already compact in $\mathcal{C}$.

\item 
	Since we have 
	$\varphi(t, K \times \Omega_\mathcal{L}) \in \tilde{K} \times \Omega_\mathcal{L}$
	so that we have `entrance' on the tail part $\Omega_\mathcal{L}$,
	then $\Omega_\mathcal{L}$ can serve as a tail in h-sets with tails,
	and the covering relations will depend only 
	on the finite-dimensional part. For small $\epsi$ we expect that
	the relation will be the same as the relation of $K$ to $\tilde{K}$ 
	in the ODE model.
\end{enumerate}

More formally, we have:
\begin{pro}
\label{pro:Omega_L}
Assume $K$ is compact. Then the set $W = K_{\mathcal{L}(L)}$ has the following properties:
\begin{enumerate}
\item 
	$W$ is compact;
\item 
	if $K$ is convex then $W$ is convex;
\item 
    let $u$ be a solution to Eq.~\eqref{eq:dde} with $f$, $g$ satisfying Condition A (see page \pageref{condA}), and let $L \ge M_f + \epsi M_g$. 
	If $\bar{u}_0 = (\bar{u}(0), \eta_0) \in K \times \Omega_{\mathcal{L}(L)} = W$ 
	then for $\bar{u}_t = (u(t), \eta_t)$ 
	we have $\eta_t \in \Omega_{\mathcal{L}(L)}$.
\end{enumerate}
\end{pro}
\noindent\textit{Proof}: 
(1) functions in $W$ are equicontinuous because they are Lipschitz. 
$W$ is bounded because for $\phi \in W$ we have: 
$\| \phi(t) -\phi(0) + \phi(0) \|_\infty \le \| \phi(0) \| + L \| t - 0 \| \le \| \phi(0) \| + L \tau$ 
and $\phi(0) \in K \subset \R^d$ compact.

(2) $v, u \in W$ and consider $w = (1-h) v + h u$. First, we have 
$w(0) \in K$ because $K$ is convex. Then, we have 
$\| w(t_2) - w(t_1) \| \le (1-h) \| v(t_2) - v(t_1)\| + h \| u(t_2) - u(t_1) \| \le L | t_2 - t_1|$.

(3) It follows directly from the boundedness assumption on $f$ and $g$ 
for DDE~\eqref{eq:dde} and the assumption on $\bar{u}_0$.
Consider $\| \bar{u}(t_2) - \bar{u}(t_1) \|$. If both $-\tau \le t_i \le 0$ 
then we use assumption on $\bar{u}_0$. If both $t_i \ge 0$ we use 
Lemma~\ref{lem:t-lipshitz}. If $t_2 < 0 < t_1$ (w.l.o.g) we have:
\begin{eqnarray*}
\| \bar{u}(t_2) - \bar{u}(t_1) \| 
	& = & \| \bar{u}_0(0) + \int_0^{t_2} F(\bar{u}(s), \bar{u}(s-\tau))ds - \bar{u}_0(t_1)\| \\
	& \le & \| \bar{u}_0(0) - \bar{u}_0(t_1)\| + \| \int_0^{t_2} F(\bar{u}(s), \bar{u}(s-\tau))ds \| \\
	& \le & L | 0 - t_1 | + (M_f + \epsi M_g) | t_2 - 0| \\
	& \le & L | t_2 - t_1 |.
\end{eqnarray*}
(in last step used the assumption $t_2 < 0 < t_1$).
It is easier to prove by picture in your mind: 
we simply glue together continuously two 
functions which both are Lipschitz with the same constant.

\qed

% SUBSECTION %%%%%%%%%%%%%%%%%%%%%%%%%%%%%%%%%%%%%%%%%%%%%%%%%%%%%%%%%%%%%%%%%%%
\subsection{Semiflows for DDEs}
Thorough this section we still assume Condition A (p.~\pageref{condA})
for $f$, $g$ in Eqs.~\eqref{eq:ode}~and~\eqref{eq:dde}
By $\Phi$ we denote the flow associated to
the ODE \eqref{eq:ode}. 
For DDE \eqref{eq:dde}, we define semiflow 
$\varphi : \R_+ \times \mathcal{C} \to  \mathcal{C}$
on the phasespace $\mathcal{C} = C([-\tau, 0], \R^d)$
by:
\begin{equation}
\varphi(t, \bar{u}_0) := \bar{u}_t
\end{equation}
for any solution $\bar{u}$ of \eqref{eq:dde}.

Under our rather strong assumptions about $f$ and $g$ 
(continuous, bounded, and Lipschitz; 
see Condition A on page \pageref{condA}) it 
is easy to check that 
\begin{pro}
\label{pro:semiflow-compact}
If $t \ge \tau$ then $\varphi(t, \cdot) : \mathcal{C} \to \mathcal{C}$ is a compact mapping.

Moreover, $\varphi(t, \cdot)|_{\R^d \times \Omega_{\mathcal{L}}}$ 
is a compact mapping (in $\mathcal{C}$) for any $t \ge 0$,
with the range $\R^d \times \Omega_{\mathcal{L}}$.
\end{pro}

Moreover, from Lemmas~\ref{lem:t-lipshitz} and \ref{lem:dde-dde-delta} we get:
\begin{cor}
If $t, t_1, t_2 \ge \tau$ then
\begin{eqnarray*}
\| \varphi(t_2, \bar{u}) - \varphi(t_1, \bar{u}) \|_\infty & \le & \left(M_f + |\epsi|M_g\right) \cdot |t_2 - t_1|, \\
\| \varphi(t, \bar{u}_2) - \varphi(t, \bar{u}_1) \|_\infty & \le & \left(1 + |\epsi| L_g \tau e^{L_f + |\epsi|L_g)t}\right) \| \bar{u}_2 - \bar{u}_1 \|_\infty.
\end{eqnarray*}
In particular, $\varphi_\epsi |_{[\tau, +\infty) \times \mathcal{C}}$ is locally Lipschitz.
\end{cor}
We will write $\varphi_\epsi$ if we want to emphasise what
value of $\epsi$ we are using. Note that for $\epsi = 0$
we have $\varphi_0(t, \bar{u}_0)(0) = \Phi(t, \bar{u}_0)$,
but $\varphi_0$ is formally a semiflow on an 
infinite-dimensional space $\mathcal{C}$.
Please note that the flow $\Phi$ is 
defined (presumably) for all times, while $\varphi$
is a semiflow, thus defined only for $t \ge 0$.

% SUBSECTION %%%%%%%%%%%%%%%%%%%%%%%%%%%%%%%%%%%%%%%%%%%%%%%%%%%%%%%%%%%%%%%%%%%
\subsection{Sections and Poincar\'e maps}
The following definitions are the same
as in \cite{wojcik-zgliczynski-dde-chaos}:

\begin{definition}
\label{def:section-ode}
The section $S$ for the flow $\Phi$ on $\R^d$ is
a $(d-1)$-dimensional smooth submanifold of $\R^d$
%homeomorphic to $\bigcup_{i=1}^{N} V_i$ with each $V_i \approx \R^{d-1}$
such that there exist $\delta > 0$ and 
a function $\alpha$ 
with the following properties:
\begin{enumerate}
\item
$\Phi|_{(-\delta, \delta) \times S}$ is a homeomorphism
with range $U = \Phi\left((-\delta, \delta) \times S\right)$;

\item
$\alpha : U \to \R$ is such that $S = \alpha^{-1}(0) \cap U$;

\item 
(transversality condition) $\nabla \alpha(u) \cdot f(u) \ge a > 0$ 
for all $u \in U$;

\item 
(existence/uniqueness condition) 
$\alpha\left(\Phi\left((-\delta, 0), S\right)\right) < 0$ and $\alpha\left(\Phi\left((0, \delta), S\right)\right) > 0$.
\end{enumerate}
\end{definition}

\begin{definition}
\label{def:poincare-map-ode}
Let $S$ be a section and let $\Sigma \subset S$ be an open set such 
that each point of $\Sigma$ 
returns to $S$ for the flow $\Phi$, i.e. for each $u \in \Sigma$ 
there exists $t_u > 0$ such that $\Phi(t_u, u) \in S$. 
Let $t_P(u) = \inf \{ t > 0 : \Phi(t, u) \in S \}$, then 
a map $P : \Sigma \to S$ defined as:
\begin{equation}
\label{eq:poincare-ode}
P(u) := \Phi(t_P(u), u)
\end{equation}
is called \emph{Poincar\'e map on $S$}.
\end{definition}

\begin{rem}
Under given assumptions the functions 
$t_P$ and $P$
are continuous. If $f$ in ODE~\eqref{eq:ode} is smooth,
then they are smooth, too.
\end{rem}

Now we would like to elaborate how the Poincar\'e maps
for ODEs~\eqref{eq:ode} persist (in an appropriate sense) 
in the DDE setting, under small delayed 
perturbation~\eqref{eq:dde} for small enough $\epsi$.
The following is inspired by Lemma~3.1 in \cite{wojcik-zgliczynski-dde-chaos},
but it has slightly different proof:
\begin{lemma}
\label{lem:poincare-dde-perturb}
Let $S$ be a section (Definiton~\ref{def:section-ode})
and let $P : S \to S$ be the Poincar\'e map for a flow $\Phi$,
with the return time function $t_P$ as in Definition~\ref{def:poincare-map-ode}.
Let $\varphi_\epsi$ be the semiflow for \eqref{eq:dde}, and assume
$f$, $g$ satisfy Condition A (p.~\pageref{condA}).

Then, for any compact $K \subset S$ 
and any $0 < \eta < \delta$ 
there exists $\epsi_{\eta,K} > 0$ such that for 
all $|\epsi| \le \epsi_{\eta,K}$, $\bar{u} \in K_\mathcal{C}$ and 
$t \in t_P(\bar{u}(0)) + [-\eta, \eta]$ 
we have:
\begin{eqnarray}
\label{eq:dde-in-U} 	\varphi_\epsi(t, \bar{u})(0) \ = \ \bar{u}(t) & \in & U, \\
\label{eq:dde-before} 	\alpha\left( \varphi_\epsi\left(t_P\left(\bar{u}(0)\right) - \eta, \bar{u}\right) \right) & < & 0, \\
\label{eq:dde-after} 	\alpha\left( \varphi_\epsi\left(t_P\left(\bar{u}(0)\right) + \eta, \bar{u}\right) \right) & > & 0, \\
\label{eq:dde-monotone} \frac{d}{dt} \alpha \left( \varphi_\epsi(t, \bar{u})(0) \right) & > & 0. 
\end{eqnarray}
In particular, for any $\bar{u} \in K_\mathcal{C}$ there exists a unique
time $t_{\bar{P}}(\bar{u}) \in t_P(\bar{u}(0)) + [-\eta, \eta]$ such that:
\begin{equation*}
\varphi_\epsi(t_{\bar{P}}(\bar{u}), \bar{u})(0) \in S.
\end{equation*}
\end{lemma}

\begin{figure}
    \centering
    \includegraphics{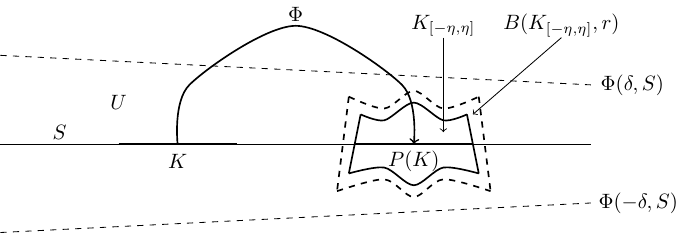}
    
    \vspace{1cm}

    \includegraphics{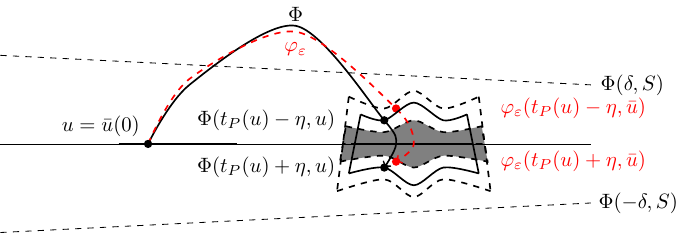}
    
    \caption{The graphical idea in the proof of Lemma~\ref{lem:poincare-dde-perturb}.
    The top picture shows the setup: the sets $K$ (left), $P(K)$ (right) on the section $S$, 
    the boundary of $U$ (dashed lines) and the sets $K_{[-\eta, \eta]}$ (solid region 
    containing $P(K)$) and $B(K_{[-\eta, \eta]}, r)$ (dashed region around $P(K)$). 
    The bottom picture shows the perturbation $\varphi_{\epsi}$ (red). For a point 
    $\bar{u} \in K_{\mathcal{C}}$ the images $\varphi(t_P(\bar{u}(0)) \pm \eta, u)$ 
    (red points) must lie close to the respective images of $P$ (black points), 
    and the $r$ is chosen in such a way, that the points must lie on the opposite 
    sides of $S$ and outside of the shaded region (grey). Together with the 
    transversality condition, it implies that $\varphi_\epsi(\cdot, \bar{u})(0)$ 
    crosses $S$ at exactly one point.}
    \label{fig:poinc-epsi}
\end{figure}

\noindent\textit{Proof}: it is similar 
to \cite{wojcik-zgliczynski-dde-chaos}, but we 
cannot alter the delay, so that the whole segment would end up in 
the strip $(-\delta, \delta) \times S$. 
Instead, we use basic estimates to show that the `head' of the solution 
$\varphi_\epsi(t, \bar{u}_0)(0)$ will be close enough to 
$\Phi(t, \bar{u}_0(0))$ for $t$ in the vicinity of $t_P(\bar{u}_0(0))$, 
for $\bar{u}_0 \in K_{\mathcal{C}}$, and
the vector field there will be just $\epsi$-small 
perturbation of the ODE system. For the graphical idea of the proof, see Figure~\ref{fig:poinc-epsi}.

We take $U = \Phi\left((-\delta, \delta) \times S\right)$ as in
Definition~\ref{def:section-ode}. Then we consider 
$$K_{[-\eta, \eta]} = \{ \Phi(t_P(u) + \eta_1, u) : u \in K, \eta_1 \in [-\eta, \eta] \} \subset U$$ 
(we think about $U$ and $K_{[-\eta, \eta]}$ to be just $S$ and 
$P(K)$ respectively `smeared out' along the flow $\Phi$). 
Obviously $K_{[-\eta, \eta]} \subset U$ and since $K_{[-\eta, \eta]}$ 
is compact and $U$ open, we know that there exists $R > 0$ such that 
for all $u \in S$ both inequalities are satisfied:
\begin{eqnarray*}
\dist\left(\Phi(t_P(u) - \eta, u), S\right) & \ge & R, \\
\dist\left(\Phi(t_P(u) + \eta, u), S\right) & \ge & R,
\end{eqnarray*}
because $\Phi$ is a homeomorphism on $U$, 
$K_{-\eta}$, $K_{\eta}$ are compact, 
$K_{\pm \eta} \cap S = \emptyset$, and $\eta < \delta$. Then, we choose an 
$0 < r < R$ so that the open set $B(K_{[-\eta, \eta]}, r) \subset U$. 

Let now $T = \max_{u \in K} \{t_P(u) + \eta\}  < \infty$ 
(as $K$ is compact and $t_P$ continuous). 
Using \eqref{eq:ode-dde-delta} we choose 
$\epsi_{K,\eta}$ such that 
\begin{equation*}
\epsi_{K,\eta} \cdot M_g \cdot T \cdot e^{L_f T} \le r.
\end{equation*}
Note that $\epsi_{K,\eta}$ depends on $\eta$ because $r$ does. 

Now, by Lemma~\ref{lem:ode-dde-delta} we have that
\begin{equation*}
\| u(t) - \bar{u}(t)\| \le r
\end{equation*}
for any $\bar{u}_0 \in K_\mathcal{C}$, where
$u$ is a solution of ODE such that $u(0) = \bar{u}_0(0)$.
Then, for any $|\epsi| \le \epsi_K$ and $t \in t_P(u) + [-\eta, \eta]$:
\begin{equation*}
\bar{u}(t) \in B(K_{[-\eta, \eta]}, r) \subset U,
\end{equation*} 
which proves \eqref{eq:dde-in-U}.

Since $r \le R$, we get $\bar{u} (t_P(\bar{u}(0)) \pm \eta) \notin S$,
hence $\alpha \circ \bar{u} (t_P(\bar{u}(0)) - \eta) < 0$
and $\alpha \circ \bar{u} (t_P(\bar{u}(0)) + \eta) > 0$,
which implies \eqref{eq:dde-before} and \eqref{eq:dde-after}.

Finally, we have 
$\alpha(\varphi(t, \bar{u}_0)(0)) = \alpha(\bar{u}(t))$
and 
\begin{eqnarray}
\frac{d}{dt} (\alpha \circ \bar{u}) (t) 
\notag			& = & \nabla \alpha (\bar{u}(t)) \cdot \bar{u}'(t) \ = \\ 
\label{eq:dt-alpha-dde}	& = & \nabla \alpha (\bar{u}(t)) \cdot \left( f(\bar{u}(t)) + \epsi g(\bar{u}(t - \tau))\right)
\end{eqnarray}
with $\| \epsi g(\bar{u}(t - \tau)) \| \le \epsi M_g$.
Since  $\grad \alpha (\bar{u}(t)) \cdot f(\bar{u}(t)) \ge a > 0$ (by assumption),
so \eqref{eq:dt-alpha-dde} is greater by $0$ by continuity, decreasing $\epsi_{K,\eta}$ 
if necessary.
\qed

\begin{cor}
Let $S$, $P$, $K$ and $\epsi_{\eta,K}$ be as in 
Lemma~\ref{lem:poincare-dde-perturb}. Then, for all $|\epsi| \le \epsi_{\eta,K}$
the first return map
$$\bar{P} : K_\mathcal{C} \ni \bar{u} \mapsto \varphi_\epsi(t_{\bar{P}}(\bar{u}), \bar{u})\  \in \mathcal{C} = C([-\tau, 0], \R^d)$$
is well defined on a section $\bar{S} = \left\{ \bar{u} \in \mathcal{C} : \bar{u}(0) \in S \right\}$.
\end{cor}
Moreover, we have:
\begin{lemma}
\label{lem:Pepsilon}
Fix $r > 0$ and let $P$, $\bar{P}$, $K$ be as in Lemma~\ref{lem:poincare-dde-perturb}.
Then, there exists $\epsi_{K,r}$, such that for all $u_0 \in K_\mathcal{C} \subset \mathcal{C}$:
\begin{equation}
\label{eq:Pepsilon}
\|P(\bar{u}_0(0)) - \bar{P}(\bar{u}_0)(0)\| \le r
\end{equation}
\end{lemma}
\noindent\textit{Proof}: 
Let $\eta < \delta$, $\epsi_{K,\eta}$ and $T$ be as in the proof 
of Lemma~\ref{lem:poincare-dde-perturb}. Moreover, assume 
$M_f \eta \le \frac{r}{2}$ and then choose
$\epsi_{K,r} < \epsi_{K,\eta}$, so that 
$\epsi_{K,r} \cdot M_g \cdot T \cdot e^{L_f T} < \frac{r}{2}$.

Fix $\bar{u}$, so that $t_P = t_P(\bar{u}_0(0))$ and 
$t_{\bar{P}} = t_{\bar{P}}(\bar{u})$ in what follows. 
Note that $t_{\bar{P}} \in t_P + [-\eta, \eta]$.
Then, for any $|\epsi| \le \epsi_{K,r}$ we have:
\begin{eqnarray*}
\|P(\bar{u}_0(0)) - \bar{P}(\bar{u}_0)(0)\| 
    & = & \| \Phi(t_P, \bar{u}(0)) - \varphi_\epsi(t_{\bar{P}}, \bar{u}) \| \\
    & \le & \| \Phi(t_P, \bar{u}(0)) - \Phi(t_{\bar{P}}, \bar{u}(0)) \| 
            + \|\Phi(t_{\bar{P}}, \bar{u}(0)) - \varphi_\epsi(t_{\bar{P}}, \bar{u}) \| \\
    & \le & M_f \eta + |\epsi|  \cdot M_g \cdot T \cdot e^{L_f T} \\            
    & \le & \frac{r}{2} + \frac{r}{2} \ = r.
\end{eqnarray*}
\qed

% SUBSECTION %%%%%%%%%%%%%%%%%%%%%%%%%%%%%%%%%%%%%%%%%%%%%%%%%%%%%%%%%%%%%%%%%%%
\subsection{Continuation of covering relations}

Now we are in a position to show that the covering relations
for ODE~\eqref{eq:ode} survive in a certain sense in DDE~\eqref{eq:dde}
for $\epsi$ small enough. 
\begin{definition}
Let $P : \mathcal{X} \to \mathcal{X}$ be a compact mapping, 
$N = \{ N_i : i = 1, \ldots, n \}$ be the finite collection of h-sets 
$N_i \subset \mathcal{X}$,
and $C \subset N \times N$. We will say that $C$ is a
\emph{collection of covering relations on $N$ for map $P$} if
for each $(N_i, N_j) \in C$ we have $N_i \cover{P} N_j$.
\end{definition}

\begin{theorem}
\label{thm:perturbation}
Let $P : S \to S$ be a Poincar\'e
map for ODE~\eqref{eq:ode},
and $C \subset N \times N$ be a collection of covering relations
on $N$ for $P$. 

For any $\epsi$ we define:
\begin{itemize}
    \item $\bar{N}^\epsi = \{ N_i \times \Omega_{\mathcal{L}(M_f + \epsi M_g)} : N_i \in N \}$;
    \item $\bar{C}^\epsi = \left\{ (\bar{N}^\epsi_i, \bar{N}^\epsi_j) : (N_i, N_j) \in C \right\}$.
\end{itemize}

Then, there exists $\epsi_{C} > 0$ such that for any $\epsi \le \epsi_{C}$,
the set $\bar{C}^\epsi$ is a collection of covering relations
on $\bar{N}^\epsi$ for the map $\bar{P}$ defined
on the section $\bar{S} \subset \mathcal{C}$, $S = \{ u \in \mathcal{C} : u(0) \in S \}$
for the semiflow of DDE~\eqref{eq:dde}.
\end{theorem}
\noindent\textit{Proof}: 
We will prove the Theorem for simplified
covering relations with one unstable direction $u = 1$,
but the same argument might be done for general covering relations and h-sets. 
The proof is a simple consequence of what was 
already said.

For a relation $N_i \cover{P} N_j \in C$ consider its homotopy $H$ 
from Definition~\ref{def:singlevalued-covering-FINITE}. For simplicity, w.l.o.g.
assume that $c_{N_i} = \operatorname{Id}$ for all $i$, so $H = H_c$.
By the definition of h-sets, 
the sets $N_j$, $N_j^+$, $H(0, N_i^-) = P(N_i^-)$ 
and $H(0, N_i) = P(N_i)$ are all compact, 
so there exists $r_{ij} > 0$ such that
$\dist(H(0, N_i^-), N_j) > r_{ij}$ and  
$\dist(H(0, N_i), N_j^+) > r_{ij}$. 

Let $K = \bigcup_i N_i$ and 
$$r_C := \min_{(N_i, N_j) \in C} \{ r_{ij} \}\text{,}$$
and choose $\epsi_C$ such that Lemma~\ref{lem:Pepsilon} holds 
with $r = r_C$ and $K$, so that $\bar{P}$ is well defined
on all of $\bar{N}^\epsi_i$ sets.

% We remind, that we can write 
% $\R^d_{\mathcal{L}} \ni u = (u(0), \omega) \in \R^d \times \Omega_{\mathcal{L}}$,
% where $\omega = u - u(0)$ and we define $\pi_\Omega(u) := u - u(0)$.

For $(\bar{N}^\epsi_i, \bar{N}^\epsi_j) \in \bar{C}^\epsi$ define a 
homotopy $\bar{H}_{ij}$ from $\bar{P}$ to $(P, 0)$ by:
\begin{equation*}
\bar{H}_{ij}(t, \bar{u}) = \left(
    (1-t)\bar{P}(\bar{u})(0) + t P(\bar{u}(0)), 
    (1-t)\pi_{\Omega} \bar{P}(\bar{u}) 
\right).
\end{equation*}
Obviously, $\bar{H}_{ij}(0, \cdot) = \bar{P}$ and $\bar{H}_{ij}(1, \cdot) = (P, 0)$. 
Since 
$$\dist\left(P(\bar{u}(0)), (1-t)\bar{P}(\bar{u})(0) + t P(\bar{u}(0))\right) 
= (1-t) \dist\left(P(\bar{u}(0)), \bar{P}(\bar{u})(0)\right) \le r < r_{ij}$$
for all $i, j$, 
we have
$\bar{H}_{ij}(t, \bar{u}) \notin \bar{N}^\epsi_j$ and 
$\bar{H}_{ij}(t, \bar{u}) \notin \bar{N}^{\epsi+}_j$ 
for $\bar{u}$ in $N_i^- \times \Omega_{\mathcal{L}}$ and
${\bar{N}^\epsi_i} = N_i \times \Omega_{\mathcal{L}}$, respectively.

From Proposition~\ref{pro:semiflow-compact} we have that
$\bar{P}$ is a compact mapping from $R^d \times \Omega_{\mathcal{L}}$ 
to $\mathcal{C}$ with range in $R^d \times \Omega_{\mathcal{L}}$, hence
$H(t, \bar{N}^\epsi_i) \subset R^d \times \Omega_{\mathcal{L}}$.
Of course all $H_{ij}(t, \cdot)$ are compact mappings, as they
are straight line homotopies between two compact maps. 

This shows properties C0-C3 of 
Definition~\ref{def:covering-with-tail} for each covering
$\bar{N}^{\epsi}_i \cover{\bar{P}} \bar{N}^{\epsi}_j$. 
To show C4 for this covering, observe that by Definition~\ref{def:singlevalued-covering-FINITE}
of the finite dimensional covering $N_i \cover{P} N_j$ 
there is a homotopy $H_{ij} : [0, 1] \times \R^d \to \R^d$ (note: no bar on $H_{ij}$),
such that $H_{ij}(0, \cdot) = P$ and $H_{ij}(1, \cdot) = (a_{ij} p, 0)$. 
Note that here we use the assumption that $c_{N_i} = \operatorname{Id}$ to simplify notation 
and to not work on $H_{c,ij}$ (see Def.~\ref{def:singlevalued-covering-FINITE} for
definition of $H_c$). Now, we 
just glue $\bar{H}_{ij}$ with the homotopy 
$(H_{ij}, 0) : [0,1] \times \R^d \times \Omega_{\mathcal{L}} \to \R^d \times \Omega_{\mathcal{L}}$ 
(transforming $(P, 0)$ to $((ap, 0), 0)$) 
in a standard way:
\begin{equation*}
H(t) = \begin{cases}
\bar{H}_{ij}(2t, \cdot) & \text{ for } t \in \left[0, \frac{1}{2}\right] \\
\left(H_{ij}(2t-1, \cdot), 0\right) & \text{ for } t \in \left[\frac{1}{2}, 1\right]
\end{cases}.
\end{equation*}
Note: $\bar{H}_{ij}(1, \cdot) =(H_{ij}(0, \cdot), 0) = (P, 0)$.

\qed

% SUBSECTION %%%%%%%%%%%%%%%%%%%%%%%%%%%%%%%%%%%%%%%%%%%%%%%%%%%%%%%%%%%%%%%%%%%
\subsection{Applications}

The perturbation Theorem~\ref{thm:perturbation} can be applied to extend any result where the covering relations
were used to prove the existence of interesting dynamics,
like chaos in Lorenz system
\cite{cap-lorenz, cap-lorenz-GiZ}, the finite dimensional 
Galerkin projection of a dissipative PDE \cite{wilczak-pde-1}, Chua's circuit 
\cite{chua-chaos-G}.

As an example, we will generalise some results of already mentioned works on
Sharkovskii property for R\"ossler ODE \cite{GZ2021, GZ2022}.
Consider the DDE~\eqref{eq:dde}
with $g= f$ and let us ask a question how small
the $\epsi$ needs to be in order to the set of coverings survives.
Of course, in such a case $f, g$ are not bounded, but we can get away 
with this assumption, if we modify the functions outside a large
enough ball around the apparent attractor. 

Consider the isolating neighbourhood $G_3$ of the attractor
for the Poincar\'e map taken from Lemma~\ref{lem:Roessler525}.
To test how small $\epsi_C$ should be, 
put $r_C = \diam(G_3) \approx 4 \cdot 10^{-4}$. This will be probably
a big overestimation, but we are interested in a rough estimate on $\epsi$.
From non-rigorous computations one can estimate very roughly that:
$T \ge 6$, $L_f = L_g \ge 10$, $M_g < 20$.
We want to have (among others):
$$|\epsi| \cdot M_g \cdot \tau \cdot e^{(L_f +|\epsi|L_g ) \cdot T} \le r.$$
For $\tau \in O(1)$, we get $|\epsi| \in O(10^{-31})$.

So, \emph{only extremely small delays don't matter}.
In the next Section we will show how to use rigorous numerics
to prove result for comparatively much bigger $\epsi$.

%% file: 4_computer_assisted.tex
%%%%%%%%%%%%%%%%%%%%%%%%%%%%%%%%%%%%%%%%%%%%%%%%%%%%%%%%%%%%%%%%%%%%%%%%%%%%%%%%
% SECTION %%%%%%%%%%%%%%%%%%%%%%%%%%%%%%%%%%%%%%%%%%%%%%%%%%%%%%%%%%%%%%%%%%%%%%
%%%%%%%%%%%%%%%%%%%%%%%%%%%%%%%%%%%%%%%%%%%%%%%%%%%%%%%%%%%%%%%%%%%%%%%%%%%%%%%%
\section{\label{sec:CAP}Computer assisted proofs}

In Section~\ref{sec:perturbation} we have shown that a sufficiently small
perturbation of an ODE with a delayed term (Eq.\ \ref{eq:dde}) 
does not destroy Sharkovskii property of the system. 
However, the maximum size of the perturbation estimated from the theory
can be extremely small, as in the exemplary R\"ossler ODE. In this section
we present an example of a computer-assisted proof
for the Sharkovskii property in a system with explicitly given parameters. 
Namely:

\begin{theorem}
\label{thm:main-numerical}
Consider the DDE \eqref{eq:dde}
where $f = g = f_{ab}$ are the R\"ossler system \eqref{eq:rossler} vector 
field with $a = 5.25$, $b=0.2$, $\tau = 0.5$, $\epsi = 0.0001$. Then
this system has $m$-periodic orbits for any $m \in \N$.
\end{theorem}
\noindent\textit{Sketch of the proof}: The proof is similar to 
\cite[Theorem~22]{GZ2022}, and specifically Lemma~\ref{lem:Roessler525} 
\cite[Lemma~21]{GZ2022}. In some Banach space $\mathcal{C}^n_p$
(defined formally in the following subsections), 
we construct explicit sets of initial segments 
$\bar{G}_3 \subset \mathcal{C}^n_p$ and $\bar{C}_i \subset 
\bar{G}_3$, $i=1,2,3$, all of them are h-sets with 
the same tail.

For the Poincare map $\bar{P} : \bar{S} \to \bar{S}$, 
$\bar{S} = \{ \phi \in \mathcal{C} : \pi_x \phi(0) = 0 \}$ 
associated to the system \eqref{eq:dde} we show with 
computer assistance that $\bar{P}(\bar{G}_3) \subsetneq \interior \bar{G}_3$ 
and $\bar{P}(\bar{C}_i) \subset \interior \bar{C}_{i+1}$. Thus we can recover 
covering relations as in Theorem~\ref{thm:sh-infD} and conclude that the system 
has periodic orbits of all periods.  

To make the article reasonably self contained and to explain how the sets
$\bar{G}_3$ and $\bar{C}_i$'s are constructed from sets $G_3$ and $C_i$ given
in Lemma~\ref{lem:Roessler525}, we need to introduce
the basic notions from \cite{szczelina-zgliczynski-focm-1, szczelina-zgliczynski-focm-2}.
The introduction is very sparse, so we highly encourage interested people 
to seek more detailed explanation in the cited papers. 

% SUBSECTION %%%%%%%%%%%%%%%%%%%%%%%%%%%%%%%%%%%%%%%%%%%%%%%%%%%%%%%%%%%%%%%%%%%
\subsection{Rigorous methods for DDEs}

The works \cite{szczelina-zgliczynski-focm-2, szczelina-zgliczynski-focm-1}
introduce a family of methods 
(denoted by $\mathcal{I}$ and $\mathcal{I}_\epsi$ for $\epsi \in \R_+$) 
to extend the initial segment of the general DDE:
\begin{equation}
\label{eq:dde-general}
u'(t) = f\left(u(t), u(t-\tau_1), \ldots, u(t-\tau_m)\right),
\end{equation}
by the method of steps, in a way similar
to how the CAPD library \cite{capd-article} 
computes enclosures on the solutions to ODEs. 
The method is not restricted to
full delays only (as in the classical method of steps), 
thus it can also be used to rigorously 
construct Poincar\'e maps on (reasonably) general sections in 
the phasespace $\mathcal{C}$ of Eq.~\eqref{eq:dde-general}. 
Of course, the perturbed system 
Eq.~\eqref{eq:dde} belongs to the class of systems \eqref{eq:dde-general}.

In short, if $X$ is a set of initial segments (given in
a certain way -- described in a moment) then we get for some $h>0$
(depending on the description of $X$) the following:
\begin{equation*}
\varphi(h, X) \subset \mathcal{I}(X)
\quad \textrm{and} \quad
\varphi(\epsi, x) \in \mathcal{I}_\epsi(X), 
\end{equation*}
for any $0 \le \epsi < h$, $x \in X$ smooth enough.
The inner workings of the algorithms are not important here,
but the description of sets is, as they are used in Theorem~\ref{thm:main-numerical}.

The set $X$ is given by the piecewise Taylor expansion 
representation of solution segments, so we need the notion of \emph{jets}. 
Let $\phi$ be some $C^{n+1}$ function,  
$\phi : [t, t + \delta) \to \R^d$ then \emph{the~forward jet of $\phi$ at $t$} is: 
\begin{equation*}
	\left(J^{[n]}_{t}\phi\right)(s) = \sum_{k=0}^{n} \phi^{[k]}(t)   s^k \quad s \in [t, t+s), 
\end{equation*}
where $\phi^{[k]}(t) = \frac{\phi^{(k)}(t)}{k!}$ are the coefficients 
of Taylor expansion ($d$ dimensional vectors). We identify the jet 
with a vector in $\R^{(n+1) \cdot d}$ 
in a standard way: 
$J^{[n]}_{t}\phi \equiv \left(\phi^{[0]}_1(t), \ldots, \phi^{[0]}_d(t) \ldots, \phi^{[n]}_d(t)\right) \in \R^{(n+1) \cdot d}$, so we can write  
$J^{[n]}_{t} \phi \in A \subset \R^{(n+1) \cdot d}$.
\begin{definition}
	The forward Taylor
	representation of $\phi : [t, t+\delta)$ $\to \R^d$ is the pair 
	$$
	(J^{[n]}_{t}\phi, \xi) \subset \R^{(n+1)\cdot d} \times C^{0}([t, t+\delta], \R^d)
	$$
	such that the Taylor formula is valid for $g$:
	\begin{equation*}
		\label{eq:forward-taylor}
		\phi(s) = \left(J^{[n]}_{t}\phi\right)(s) + (n+1) \int_t^\delta \xi(s)   (t-s)^n ds,\quad s \in [t, t + \delta).
	\end{equation*}   
\end{definition}
Obviously, $\xi(s) = \phi^{[n+1]}(s)$, and it is assumed $\|\phi^{[n+1]}(s)\|$ 
is bounded over $[t, t+\delta)$. 
Conversely, if we have any $(j, \xi) \in \R^{(n+1) \cdot d} \times \mathcal{C}([0, h], \R^d)$ 
then we say it is equivalent to some $\phi : [t, t + \delta]$, iff $J^{[n]} = j$. 
We will write $(j, \xi)(t)$ to denote the evaluation of such $\phi(t)$ 
by \eqref{eq:forward-taylor}, and we will denote the equivalence by $(j, \xi) \equiv \phi$. The quantities $t$ and $\delta$ in $\phi$ will be usually known from the context. 
\begin{definition}
	\label{def:cetap}
	For a vector $\eta = (n_1, \ldots, n_p)$ 
	the space $C^\eta_p$ is defined as:
	\begin{eqnarray*}
		C^\eta_p & = & \left\{ u : [-1, 0] \to \R: \exists j_i \in \R^{n_i + 1},  \exists \xi_i \in C^{0}([0, h], \R) \right. \\ 
		& & \left. \textrm{ such that }  u_{[t_i, t_{i-1})} \equiv (j_i, \xi_i)\right\},
	\end{eqnarray*}
	where $t_i = -i   h_p$ for $h_p = \frac{\tau}{p}$ form \emph{a grid of size $p$} over $[-\tau, 0]$. 
\end{definition}
We will usually drop subscript $p$ in $h$, as $p$ is fixed in a given context.
We will write $C^{n}_{p}$ if $\eta = (n, \ldots, n)$.

\begin{definition}
	A ($p$, $\eta$)-representation (or just the representation)
	of a function $u \in C^\eta_p$ is the collection 
	of $(z, j_1, \ldots, j_p, \xi_1, \ldots, \xi_p)$,
	where $z = u(0)$ and $u_{[t_i, t_{i-1}]} \equiv (j_i, \xi_i)$
	as in Definition~\ref{def:cetap}.
\end{definition}
The vector $(z, j_1, \ldots, j_p, \xi_1, \ldots, \xi_p) = (z, j, \xi)$
can be thought of as an element of 
$\R^M \times (C^0([0, h], \R^d))^p = \R^M \times (C^0([0, h], \R))^{d\cdot p}$,
with $M = M(p, \eta) = d \cdot (1 + \sum_{j=1}^{p} (\eta_i + 1))$.
This also shows that $C^\eta_p$ is a Banach space, with the
norm given by the maximum norms in $\R^M$ and $C^0([0, h], \R)$) 
in a standard way: 
$$\| (z, j, \xi) \|_{C^\eta_p} = \| (z, j) \|_1 + \sum_{j=0}^{d \cdot p} \| \xi_j \|_\infty.$$ This will be used in Section~\ref{sec:CAP}.

Finally, let $Z, Y$ be sets in $\IR^d$, i.e.
the products of closed intervals.
For $\xi \in C^0(Z, \R^d)$ we write $\xi \in Y$ to denote the fact that:
$$\left[\inf_{t \in Z} \xi_j(t), \sup_{t\in Z} \xi_j(t)\right] \subset Y_j, \quad \forall j = 1,\ldots,d.$$ 
\begin{definition}
	\emph{The representable ($d$, $p$, $\eta$)-functions-set} (or just \emph{f-set})
	is a pair $(A, \Xi) \subset \R^M \times \IR^{p \cdot d}.$ The support 
	of the f-set is $\mathfrak{X}(A, \Xi) \subset C^\eta_p$: 
	\begin{equation}
		u \in \mathfrak{X}(A, \Xi) \iff \exists (z, j) \in A, \exists \xi \in \Xi \textrm{ such that } u = (z, j, \xi).
	\end{equation}
\end{definition}
We will often identify an f-set with its support. For $X = \mathfrak{X}(A, \Xi)$ 
we will use projections: $\mathfrak{z}(X)$, $\mathfrak{j}(X)$, 
$\mathfrak{j}_i(X)$, $\xi(X)$, $A(X)$. 
The projection $\mathfrak{z}(X)$ is called \emph{head} and it is
a $d$-dimensional vector that corresponds to the value $u_t(0)$
of a solution to the DDE~\eqref{eq:dde}, and 
$\pi_z(\mathfrak{z}(X))$ is its third coordinate in
the \emph{ambient} $d$-dimensional space of the ODE state space. 

For $u \in C^\eta_p$ by $\mathfrak{X}(u) := \mathfrak{X}(\{a\}, \Xi)$ 
we will denote the smallest f-set that contains the representation 
of $u \equiv (a, \xi)$, $\xi \in \Xi$. 

For any set $\mathfrak{X}(A, \Xi)$ the space $C^\eta_p$ will
be always known from the context as the algorithms 
$\mathcal{I}$, $\mathcal{I}_\epsi$ and Poincar\'e
maps $P : C^n_p \to C^n_p$ in the \texttt{capdDDEs} 
library presented in \cite{szczelina-zgliczynski-focm-2} take care
of tracking the $\eta$'s and making sure that all discontinuities are 
harmless for the computations to be rigorous. In case the computations 
cannot be make rigorous, the library raises exceptions and 
the computation results in a failure. The methods are built upon 
the foundation of the CAPD library \cite{capd-article}
that supports interval arithmetic, rigorous algebraic methods, etc. 

On the other hand, in the context of covering relations in infinite 
dimensions, the sets $\mathfrak{X}(A, \Xi)$
will be the h-sets with tails, where naturally $A \subset \R^M =: \mathcal{X}_1$ 
will be the h-set in a finite dimensional space $\R^M$, whereas 
the tail will be $\Xi \subset \left(C([0, h], \R)\right)^{d \cdot p}$
with the maximum product norm of supremum norms in each of $C([0, h], \R)$.

% SUBSECTION %%%%%%%%%%%%%%%%%%%%%%%%%%%%%%%%%%%%%%%%%%%%%%%%%%%%%%%%%%%%%%%%%%%%%%%

\subsection{\label{sec:preparing}Preparing initial data}

This section is technical and is included
for people interested in the details of 
computer assisted proofs \cite{dde-ros-sha-codes}, or in
reproducing the results of this study. In this section 
we present the way how we generated the data for our proof, as it is not feasible to supply the initial
conditions by hand. 

Assume $p, n \in \N$ are fixed parameters of $C^n_p$ space.
We remind the reader that a function $u \in C^n_p$ is represented
with a triple $(z, j, \xi)$, where $u(0) = z \in \R^d$, $j$ describes the 
jets of the solution segment at grid points and constitutes the 
representation of $u$ on the interval $[-\tau, 0)$ (open from the right), 
and $\xi$ represents estimates on $(n+1)$st derivative of this 
solution segment. The pair $(u(0), j)$ can be
seen as a vector in $\R^M$ with $M = M(d, p, n)= d \cdot (p \cdot (n+1) + 1)$,
so in the sequel we may apply the standard projections to elements $u \in C^n_p$:
$\pi_i u := \pi_i (u(0), j)$. For $i = 1, \ldots, d$ we have of course 
$\pi_i u = \pi_i u(0)$. 

The dimension of the ambient space in the case of
our toy example of R\"ossler equation and its DDE perturbation
is $d = 3$. In the computer assisted proof we will use 
$n=3$ and $p = 32$, and thus $M = 387$.

We aim to define our h-sets with tails in the $C^n_p$ space using the
finite dimensional sets defined in Lemma~\ref{lem:Roessler525}. We expect that 
the heads $\mathfrak{z}(u)$ of the solutions for the DDE Poincare map $\bar{P}_\epsi$
will not diverge too much from the images of the ODE map $P$, so
we expect the covering relations on the heads $\bar{P}(\bar{u})(0)$ 
to be preserved in a manner similar to what was proven
analytically for arbitrarily small $\epsi$
in Section~\ref{sec:perturbation}.

We need to consider how to define the tail (that is, the representation of $u$ on $[-\tau, 0)$) of our sets (it
will be the same for all sets $\bar{G}_3$ and $\bar{C}_{1,2,3}$).
Since the dimension $M$ of the representation is a big number,
it will not be feasible to present data for proofs in the 
manuscript and not even in the code of computer programs. It is also not possible
to provide data by hand, therefore in the proof they are stored in 
external files (see the directory \texttt{proof-data} in \cite{dde-ros-sha-codes}).
Here, we will describe a heuristic procedure
to obtain the sets $\bar{G}_3$ and $\bar{C}'_i$, $i = 1, 2, 3$ from
respective sets given in Lemma~\ref{lem:Roessler525}. 

\begin{figure}
	{\begin{center}
	    \includegraphics[width=0.6\linewidth]{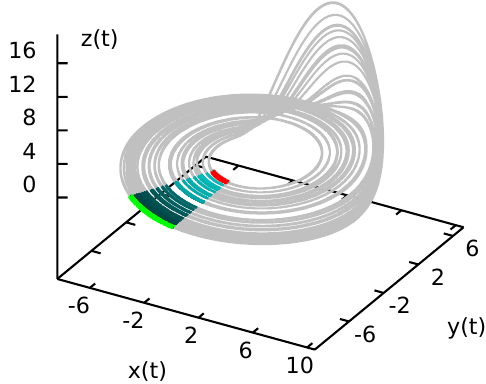}
	\end{center}} \vspace*{-1.5em} 
	\caption{\label{fig:attractor}The figure shows the long trajectory to the R\"ossler system for the initial point $(0,-5, 0.03)$ (gray). Note that this is not the attractor which is shown in Figure~\ref{fig:attr3}.
	The trajectory is drawn in the ambient space $\R^3$, together with the approximation of the set $\bar{\mathcal{A}}$
	of segments on the section $\bar{S}$ (coloured in a blue gradient). There are 100 segments used to estimate the \emph{tail part} in the heuristic procedure proposed in Section~\ref{sec:preparing}. Each segment spans the time of the full delay $\tau = \frac{1}{2}$. The red and green segments are $r$ and $l$, respectively, that are used in the heuristic procedure to define vector $\bar{m}_1$. }
\end{figure}

To get some intuition, consider the Poincar\'e map 
$P$ for ODE~\eqref{eq:ode} and assume there is a closed and 
bounded set $G \subset S \subset \R^d$ for the map $P$, 
such that $P|_G$ is a homeomorphism on its image $P(G) \subset G$
(as in our toy example with $G = G_3$). 
Let $\mathcal{A} = \operatorname{Inv}(G, P)$
be the invariant part of $G$.
Now, consider DDE~\eqref{eq:dde} with $\epsi = 0$. 
It is easy to see that the solutions to both DDE and ODE systems are the same
(treated as functions $\R \to \R^d$) and we have the following:
\begin{pro}
\label{pro:barmathcalA}
    The set 
    \begin{equation}
    \label{eq:Attractor}
    \bar{\mathcal{A}} = \left\{ \Phi(\cdot,u)|_{[-\tau, 0]} : u \in \mathcal{A} \right\}\ \subset \mathcal{C},
    \end{equation}
    is a closed subset of a convex and compact (in the phasespace $\mathcal{C}$) 
    set $N = G \times \Omega_{\mathcal{L}(M_f)}$. 
    Moreover, for any 
    $u = (u(0), \zeta) \in N$
    with $u(0) \in \mathcal{A}$, we have $\bar{P}_0(u) \in \bar{\mathcal{A}}$,
    no matter what $\zeta$ is. 
    
    In other words, $\bar{\mathcal{A}}$ and $N$ are 
    forward invariant under $\bar{P}_0 : N \to N \subset \mathcal{C}$ and
    $N$ is the trapping region of the trapped attractor 
    $\bar{\mathcal{A}}$ \cite{attractor-scholarpedia}.
\end{pro}
You can see an approximation of the set $\bar{\mathcal{A}}$ for
our toy example in Figure~\ref{fig:attractor}. 

In computations we will we need to consider the 
dynamics of DDE~\eqref{eq:dde} and $\bar{P}_\epsi$ 
in the space $C^n_p$.
Assume now that $f$ and $g$ in the DDE~\eqref{eq:dde} 
are $C^\infty$ functions and for $\epsi \in [0, \epsi_0]$ 
the perturbed map $\bar{P}_\epsi$ is such
that $\bar{P}_\epsi(C^n_p) \subset C^n_p$. The easiest way
to guarantee that is 
just to require that the return time $t_{\bar{P}_{\epsi}} > n \cdot \tau$,
i.e. it is \emph{long enough} for the solutions to be smoothed
by the DDE dynamics \cite{szczelina-zgliczynski-focm-2}. 
Please note that $\bar{P}_0(C^n_p) \subset C^n_p$
always since the solutions to ODE~\eqref{eq:ode}
are smooth, so obviously we have:
$\bar{\mathcal{A}} \subset C^n_p$.  

We expect that the set $\bar{\mathcal{A}}$ survives
for $\bar{P}_{\epsi}$, with $\epsi$ small enough, or more 
precisely, we expect that its trapping regions defined
in $C^n_p$ survive:
\begin{pro}
    Assume $f, g$ are smooth in DDE~\eqref{eq:dde} and 
    there are given $\epsi_0 > 0$ and a bounded 
    nonempty set $\bar{G} \subset C^n_p$ 
    such that $\bar{\mathcal{A}} \subset \interior \bar{G}$.
    
    Then, there exists $0 < \epsi_1 \le \epsi_0$ such that
    for all $|\epsi| \le \epsi_1$ we have:
    $\bar{P}_\epsi(\bar{G}) \subset \interior \bar{G}$.
\end{pro}
\noindent\textit{Proof}: 
it is an easy consequence of Proposition~\ref{pro:Omega_L}
and Lemma~\ref{lem:Pepsilon} (with $K = G$ and $r = \dist(\bar{\mathcal{A}}, \bd \bar{G})$ so that $\epsi_1 = \epsi_{K,r}$).
\qed

In practice, we will fix the desired $\epsi$ and 
construct 
some set $\bar{G}$ such that $\bar{P}_\epsi (\bar{G})\subset \bar{G}$. 
We will show now how to find a reasonable set $\bar{G}$ for our toy example. The procedure is general and might be applied for other systems as well.
We remind that:
\[
G_3= \bmatrix -6.38401 \\ 0.0327544 \endbmatrix
+ \bmatrix -1.& 0.000656767 \\ -0.000656767 & -1. \endbmatrix \cdot
\bmatrix \pm 3.63687 \\ \pm 0.0004 \endbmatrix
\ =: \ 
u_0 + M \cdot r_0
\text{.}
\] 
We will utilise the fact that the invariant set  
$\mathcal{A}$ is almost one dimensional and in consequence the set $\bar{\mathcal{A}}$ is close to being a hyperplane in the space $\mathcal{C}$. 
We also would like to make $\bar{G}_3$ such that $\mathfrak{z}(\bar{G}_3) = G_3$,
so let $M = \bmatrix m_1 & m_2 \endbmatrix$ with 
$m_i$ being the $i$th column of $M$. We will use them in a moment.
In what follows we drop the subscript in $G_3$, $\bar{G}_3$ 
for convenience, and we will write just $G$ and $\bar{G}$. 

\begin{rem}
One might be tempted to utilise Eq.~\eqref{eq:Attractor}
to define $\bar{G}$ as the following set:
\[
W = \left\{ \Phi(\cdot,u)|_{[-\tau, 0]} : u \in G \right\},
\]
This, however, is not a great idea, as the flow $\Phi$ in our case
is strongly expanding when going backward in time, making the set $W$
extremely large compared to its expected 
image $\bar{P}_0(W) \approx \bar{\mathcal{A}}$.
\end{rem}

Our heuristic procedure would be to find 
vectors $\bar{u}_0, \bar{m}_1, \bar{m}_2 \in C^n_p$
such that $r_* := \dist(\bar{\mathcal{A}}, \bar{u}_0 + \spans(\bar{m}_1))$
is minimal. Assume for now that:
\begin{enumerate}
    \item $\pi_{y,z} (\mathfrak{z}(\bar{u}_0)) = u_0$,
    \item $\pi_{y,z} (\mathfrak{z}(v)) = m_1$,
    \item $\pi_{y,z} (\mathfrak{z}(\bar{m}_2)) = m_2$ and $m_2(t) = 0$ for $t \in [-\tau, 0)$.
\end{enumerate}
Let 
\[ 
    \bar{G} = \bar{u}_0 + \bmatrix \bar{m}_1 & \bar{m}_2 \endbmatrix \cdot r_0 + \cB_{\spans(\bar{m}_1,\bar{m}_2)^\bot}(r_*).
\]
We see that $\mathfrak{z}(\bar{G}) = G$. If conditions
1 and 2 are not satisfied, one can modify
$\bar{u}_0$ and $\bar{m}_1$ slightly, maybe 
at the expense of increasing $r_*$ a little bit. 

Such a set $\bar{G}$ would be better than 
simply following Section~\ref{sec:perturbation}
and choosing $\bar{G}' \in C^n_p$ as 
$$\bar{G}' = \mathfrak{X}\left(G \times \Pi_{i=3}^{M} [l_i, r_i] \times \Pi_{j=0}^{d \cdot p} [L_i, R_i]\right),$$ 
with $l_i, r_i, L_j, R_j \in \R$ big enough, so that $\mathfrak{X}(\bar{u}) \in \bar{G}$
for all $\bar{u} \in \bar{\mathcal{A}}$. It is easy to see,
as the set $\bar{G}'$ is basically the set $\bar{G}$ defined
with the vector $m_1(t) \equiv 0$ on $[-\tau, 0)$, so it must be worse
than the choice with optimal $r_*$. In practice it would be a lot worse. 

The program \texttt{Roessler\_DDE\_gencoords.cpp} from \cite{dde-ros-sha-codes}
selects the vector $\bar{m}_1$ in the following way (see Figure~\ref{fig:attractor}
to get some intuition):
\begin{enumerate}
    \item A long trajectory $\bar{u}$ is computed non-rigorously, so that the set
    $\hat{\mathcal{A}} = \{ \bar{u}_t : \bar{u}_t(0) \in S\} \subset \bar{S}$ 
    consists of at least 100 loosely distributed points. The set $\hat{\mathcal{A}}$ should be a 
    reasonable approximation of $\bar{\mathcal{A}}$. The initial point of the 
    trajectory is $\bar{u}(0) = (0, -5, 0.03)$.
    
    \item The vectors $\bar{l}, \bar{r}$ from $\hat{\mathcal{A}}$ are selected 
    so that $\pi_y \mathfrak{z}(\bar{l})$ is minimal and $\pi_y \mathfrak{z}(\bar{r})$ 
    is maximal, then the vector $v'$ is computed as:
    $$v' = \frac{\bar{r} - \bar{l}}{\pi_y z\left(\bar{r} - \bar{l}\right)}.$$
    This vector should be a good approximation of the optimal $\bar{m}_1$.
    
    \item The vector $m_1$ (\texttt{ydir} in the program) is defined as 
    \begin{eqnarray}
    \label{eq:choice-m1}
    && \pi_1 \bar{m}_1 = 0, \quad \pi_{2,2} \bar{m}_1 = m_1, \\
    \notag
    && \pi_{j > 3} \bar{m}_1 = \pi_{j > 3} v'.
    \end{eqnarray}
\end{enumerate}
Similarly, the vector $\bar{u}_0$ is selected as the vector 
from the trajectory with $\mathfrak{z}(\bar{u}_0)$ closest to $u_0$ and its
head is modified to match $u_0$. Together with the choices of
$m_1$ in Eq.~\eqref{eq:choice-m1} and the definition of $m_2$ 
those modifications guarantee that $\pi_{y, z}(\bar{G}) = G$.

In order to construct the sets $\bar{C}'_i$, 
the relative positions of middle points of sets $\bar{C}'_i$
with respect to hyperplane $\bar{u}_0 + \spans(\bar{m}_1)$ are computed
so that their heads again coincide with the middle points of sets $C'_i$.

Next, the program \texttt{Roessler\_DDE\_prepare\_tail.cpp} 
computes $r_*$ (or more precisely, tight enclosures on $\bar{\mathcal{A}}$
projected onto $\bar{u}_0 + \spans(\bar{m}_1)$). The computations
are done in parallel, mimicking the subdivisions made by the final
program \texttt{Roessler\_DDE\_proof\_piece.cpp} used in the proof
of Theorem~\ref{thm:main-numerical}.

% SUBSECTION %%%%%%%%%%%%%%%%%%%%%%%%%%%%%%%%%%%%%%%%%%%%%%%%%%%%%%%%%%%%%%%%%%%
\subsection{Proof of Theorem~\ref{thm:main-numerical}}

In what follows, the $\bar{X}$ for a given set $X$ 
is obtained as described in the previous subsection.
The proof is computer assisted with source codes available in 
\cite{dde-ros-sha-codes}. 

To represent functions in the phasespace $C^n_p([\tau, 0], \R^d)$
we choose $p = 32$, $n = 3$, $d = 3$ and 
the delay $\tau = \frac{1}{2}$.

The computer program verifies the following inclusions:
\begin{itemize}
    \item $\bar{P}(\bar{G}_3) \subset \bar{G}_3$,
    \item $\bar{P}(\bar{C}'_1) \subset \bar{C}'_2$,
    \item $\bar{P}(\bar{C}'_2) \subset \bar{C}'_3$,
    \item $\bar{P}(\bar{C}'_3) \subset \bar{C}'_1$,
\end{itemize}
under the first return map $\bar{P} : \bar{S} \to \bar{S} \subset C^3_{32}\left(\left[-\frac{1}{2}, 0\right], \R^3\right)$, with   
$$\bar{S} = \left\{ u \in C^3_{32}\left(\left[-\frac{1}{2}, 0\right], \R^3\right): \pi_x(u(0)) = 0 \right\}.$$

\begin{figure}
	\includegraphics[width=0.8\linewidth]{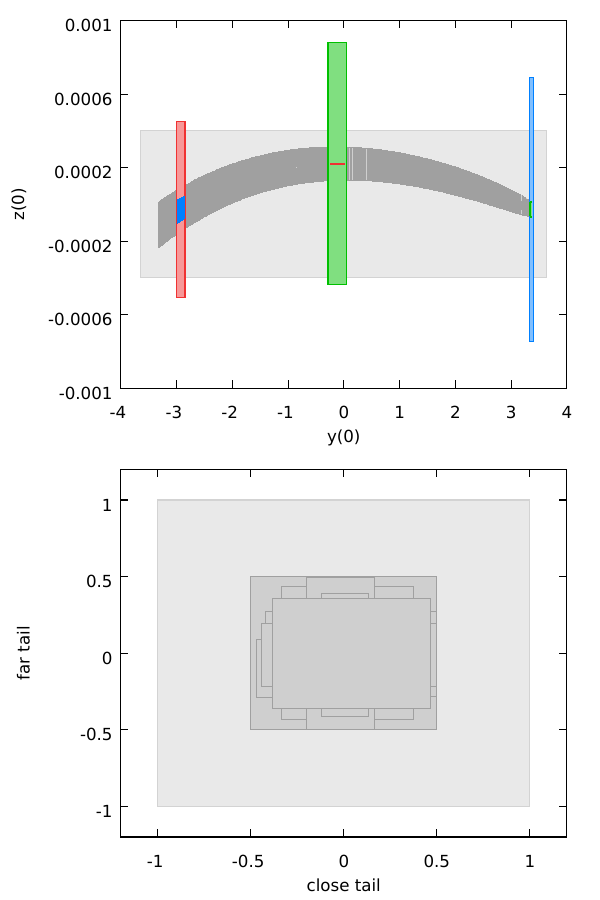}
	\caption{\label{fig:cap}The graphical representation of the computer
	assisted proof of Theorem~\ref{thm:main-numerical}. The sets are presented 
	in the good coordinate frame given by $M$ of set $G_3$, and therefore they are rectangles (cf. Fig~\ref{fig:grid3}, where the similar sets are presented in the original $(y, z)$ coordinates). On top, the projection onto $(y, z)$ coordinates
	of the heads is shown for all the sets, whereas in the bottom picture, the projections onto the close and far tails in the max-norm are presented. 
	The sets $\bar{G}_3$, $\bar{C}'_i$, $i=1,2,3$ are
	plotted in grey, red, blue and green, respectively (cf. Fig~\ref{fig:grid3}).
	The images of the sets are plotted in a more saturated respective colour.}
\end{figure}

The sets are defined as in Section~\ref{sec:preparing}.
The rough graphical representation of the verification
procedure is presented in Figure~\ref{fig:cap}. 

Now, if we define $C_i = C'_i \cap G_3$,
then $C_i \subset C_{i+1}$, and $C_i \subset G_3 \subset \domain(P)$ 
fulfil the assumptions of 
Theorem~\ref{thm:sh-infD}. Thus, the system \eqref{eq:dde} with the
specified values of parameters has a 3-period orbit and, in turn,
all $m$-periodic orbits for $m \in \N$.
\qed

As a last remark, we would like to point out that the value of $\epsi$ 
from Theorem~\ref{thm:main-numerical} is of 25 orders of magnitude larger
than the theoretical upper bound estimated in section \ref{sec:perturbation}.

%% file: 5_conclusions.tex
%%%%%%%%%%%%%%%%%%%%%%%%%%%%%%%%%%%%%%%%%%%%%%%%%%%%%%%%%%%%%%%%%%%%%%%%%%%%%%%%
% SECTION %%%%%%%%%%%%%%%%%%%%%%%%%%%%%%%%%%%%%%%%%%%%%%%%%%%%%%%%%%%%%%%%%%%%%%
%%%%%%%%%%%%%%%%%%%%%%%%%%%%%%%%%%%%%%%%%%%%%%%%%%%%%%%%%%%%%%%%%%%%%%%%%%%%%%%%
\section{\label{sec:conclusions}Final remarks}

\begin{figure}[!htbp]
	{\begin{center}
	    \includegraphics[width=0.6\linewidth]{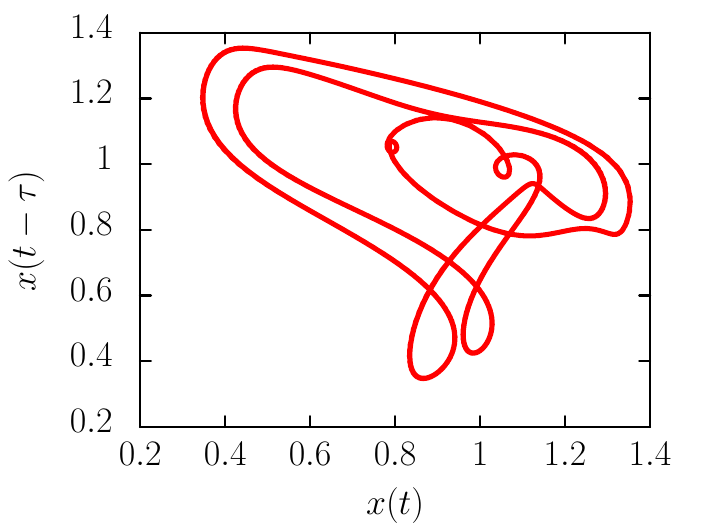}
	\end{center}} \vspace*{-1.5em} \caption{\label{fig:MG}Projection of the expected attracting $3$-periodic orbit for the system \eqref{eq:MG} with $\gamma=1$, $\beta = 2$, $\tau = 2$ and $n=9.7056$.}
\end{figure}

The Sharkovskii property for a dynamical system with a 3-periodic orbit 
is a chaotic-like behaviour of the 
system in the sense of Li-Yorke chaos \cite{LiYorke} and therefore our method may 
be seen as a step forward to proving symbolic dynamics \cite{Morse}. The existence 
of symbolic dynamics for  R\"ossler system with `classical' values of parameters 
$a=5.7$, $b=0.2$ was proven with assistance of \texttt{CAPD} library 
in \cite{ZRossler} and for a DDE-perturbed system in \cite{szczelina-zgliczynski-focm-2}.
The more challenging problem would be to prove the Sharkovskii property 
for the well known Mackey-Glass equation \cite{mackey-glass}: 
\begin{equation}\label{eq:MG}
x'(t)= \dfrac{\beta x(t-\tau)}{1+(x(t-\tau))^n}-\gamma x(t)
\end{equation}
which is expected to have an attracting 3-periodic orbit for 
$\gamma=1$, $\beta = 2$, $\tau = 2$ and $n=9.7056$ (see Fig.~\ref{fig:MG}) 
and numerical studies suggest that for those values of parameters
the apparent attractor is indeed one dimensional \cite{Humphries-Duruisseaux-mg-numerical-biff}.
In this case, however, the choice of a convenient Poincar\'e section
and the attractor's isolating neighbourhood seems to be extremely difficult 
and lies beyond the scope of this paper. But 
we would like to stress the fact that the Theorem~\ref{thm:sh-infD}
is already in a good shape to be applied in \eqref{eq:MG}. 

